\documentclass[12pt, reqno]{amsart}

\usepackage{amsmath, amssymb, amsthm, amsfonts, mathtools, array, xfrac, float, multirow, indentfirst, url, enumerate, comment, fullpage, afterpage, pifont, makecell}

\mathtoolsset{showonlyrefs}

\usepackage{xcolor}
\usepackage{hyperref}
\usepackage{diagbox}
\usepackage{booktabs, siunitx, xcolor, graphicx}
\usepackage{amsmath}
\usepackage{array}
\newcolumntype{?}{!{\vrule width 1pt}}
\usepackage{appendix}
\usepackage{tikz}
\usepackage[outline]{contour} % glow around text
\usetikzlibrary{decorations.markings}

\tikzset{>=latex} % for LaTeX arrow head
\colorlet{myblue}{blue!80!black}

\allowdisplaybreaks

\theoremstyle{plain}
\newtheorem{theorem}{Theorem}[section]

\newtheorem{lemma}[theorem]{Lemma}

\newtheorem{proposition}[theorem]{Proposition}

\theoremstyle{definition}
\newtheorem*{remark}{Remark}

\renewcommand{\Re}{\operatorname{Re}}
\renewcommand{\Im}{\operatorname{Im}}

\title[Mertens function and the reciprocal of the Riemann zeta function]{New explicit bounds for Mertens function and the reciprocal of the Riemann zeta function}

\author[E.~S.~Lee]{Ethan~Simpson~Lee}
\address{ESL: University of the West of England, School of Computing and Creative Technologies, Coldharbour Lane, Bristol, BS16 1QY} 
\email{ethan.lee@uwe.ac.uk}
\urladdr{\url{https://sites.google.com/view/ethansleemath/home}}

\author[N.~Leong]{Nicol~Leong}
\address{NL: University of Lethbridge, Alberta, Canada} 
\email{encinicol.leong@uleth.ca}

\begin{document}

\begin{abstract}
We prove new unconditional and explicit upper bounds for the reciprocal of the Riemann zeta function in the critical strip, of the orders
\begin{equation*}
    (\log t)^{11/12}, \quad
    (\log t)^{11/12}(\log\log t)^{-3/4}, \quad\text{and}\quad
    (\log t)^{2/3}(\log\log t)^{1/4} ;
\end{equation*} 
these hold in prescribed classical, Littlewood, and Korobov--Vinogradov zero-free regions, respectively. From these bounds, we also derive new unconditional and explicit upper bounds for the Mertens function $M(x)$, of the orders
\begin{align*}
    x (\log x)\exp\!\left(-\eta_1 \sqrt{\log x}\right), \qquad x (\log x) \exp\!\left(-\eta_2 \sqrt{\log x \log\log x}\right),
\end{align*}
and
\begin{equation*}
    x (\log x)\exp\!\left(-\eta_3 (\log x)^{3/5}(\log\log x)^{-1/5}\right),
\end{equation*}
 for suitable constants $\eta_1,\eta_2,\eta_3>0$. A key feature of our approach is the use of computational smoothing, which ensures the resulting numerical bounds are strong.
\end{abstract}

\maketitle

\section{Introduction}

Recall that the Riemann Hypothesis (RH) postulates every non-trivial zero $\varrho = \beta + i\gamma$ of the Riemann zeta function $\zeta(\sigma + it)$ satisfies $\beta=1/2$ and $\gamma\neq 0$. The Riemann height is a constant $\mathcal{H} > 0$ such that the RH is verified to be true for every $\varrho$ satisfying $|\gamma|\leq \mathcal{H}$; Platt and Trudgian \cite{PlattTrudgianRH} have verified that
\begin{equation}\label{eqn:PlattTrudgian}
    \mathcal{H} = 3\,000\,175\,332\,800 
\end{equation}
is admissible. In the region $| t | > \mathcal{H}$, where the RH is not verified, it is also known that $\zeta(\sigma + it)\neq 0$ in the three following forms of zero-free region:
\begin{equation}\label{zfr regions}
    \sigma \ge 1- \frac{1}{Z_1\log t},\quad 
    \sigma \ge 1- \frac{\log\log t}{Z_2 \log t},\quad\text{and}\quad 
    \sigma\ge 1- \frac{1}{Z_3(\log t)^{2/3}(\log\log t)^{1/3}} ,
\end{equation}
for some positive constants $Z_1,Z_2,Z_3$. Respectively, these regions are called the {classical}, {Littlewood}, and {Korobov--Vinogradov} zero-free regions. If $|t|\ge 3$, then we know $\zeta(\sigma + it)\neq 0$ in the regions defined in \eqref{zfr regions} with 
\begin{equation}\label{z zerofree const}
    Z_1 = 5.558691,\quad 
    Z_2 = 21.233,\quad\text{and}\quad 
    Z_3 = 53.989.
\end{equation}
For more details about the computations in \eqref{z zerofree const}, see \cite{mossinghoff2024explicit, yang2024explicit, bellotti2024explicit}, respectively. 

A big-$O$ bound is called effective if all implied constants are computable, and explicit if these constants are concretely described. 
In lieu of the RH, upper bounds for the reciprocal of the Riemann zeta function, $\zeta(\sigma + it)^{-1}$, can only be obtained inside zero-free regions like those in \eqref{zfr regions}. Such bounds have numerous applications, as $\zeta(\sigma + it)^{-1}$ often appears naturally in Perron-type formulae and related contour-integral arguments, where upper bounds for $|\zeta(\sigma + it)|^{-1}$ are required to move the line of integration. In addition, it is well known that $M(x) = o(x)$ is equivalent to the prime number theorem and that the RH is equivalent to $M(x) \ll x^{1/2 + \varepsilon}$ for any $0<\varepsilon < 1/2$, so bounds for $M(x)$ are deeply connected to the distribution of primes and zeros of $\zeta(s)$.

Notably, from explicit upper bounds for $|\zeta(\sigma + it)|^{-1}$ within zero-free regions of the forms given in \eqref{zfr regions}, one can establish explicit bounds for $M(x)$, defined as
\begin{equation}\label{eqn:Mx_def}
    M(x) = \sum_{n \le x}\mu(n) ,
\end{equation}
where $\mu(n)$ denotes the M\"{o}bius function. For suitable constants $\eta_1,\eta_2,\eta_3>0$, these bounds can take the asymptotic forms
\begin{equation}\label{eqn:explicitise_me}
\begin{split}
    M(x) &\ll x (\log x)\exp\!\left(-\eta_1 \sqrt{\log x}\right),\\
    M(x) &\ll x (\log x) \exp\!\left(-\eta_2 \sqrt{\log x \log\log x}\right), \quad\text{and}\\
    M(x) &\ll x (\log x)\exp\!\left(-\eta_3 (\log x)^{3/5}(\log\log x)^{-1/5}\right),
\end{split}
\end{equation}
corresponding respectively to zero-free regions of the forms described in \eqref{zfr regions}. 

In this paper, we prove new unconditional and explicit upper bounds for $\zeta(s)^{-1}$ of orders
\begin{equation}\label{eqn:asymptotic_forms}
    (\log{t})^{11/12} , \quad
    (\log{t})^{11/12} (\log\log{t})^{-3/4} , \quad\text{and}\quad
    (\log{t})^{2/3} (\log\log{t})^{1/4} ,
\end{equation}
inside prescribed classical, Littlewood, and Korobov--Vinogradov zero-free regions. Many of these results provide the first explicit bounds of their kind. Moreover, our computational approach yields numerically strong results that improve upon the latest explicit bounds established in \cite{CullyHugillLeong,nicol,Trudders}, both numerically and asymptotically. Subsequently, we prove new unconditional and explicit upper bounds for $M(x)$ of the forms in \eqref{eqn:explicitise_me}. 

Since our full results are somewhat technical to state, we begin by presenting a headline theorem. This theorem is a special case of the more extensive computations presented in Theorems \ref{thm:objective1} and \ref{thm:objective2}, in Section~\ref{sec:results}, which also contains the associated technical details. 

\begin{theorem}\label{thm:headline}
Let $t \geq \mathcal{H}$ with $\mathcal{H}$ defined as in \eqref{eqn:PlattTrudgian}. 
First, if 
\begin{equation}\label{eqn:classical-obj1}
    \sigma \ge 1-\frac{1}{5.56 \log{t}} ,
    \quad\text{then}\quad
    \left|\frac{1}{\zeta(\sigma+it)}\right| \le 1\,237 (\log{t})^{11/12} .
\end{equation}
Second, if 
\begin{equation}\label{eqn:Littlewood-obj1}
    \sigma \ge 1-\frac{\log\log{t}}{21.24 \log{t}} ,
    \quad\text{then}\quad
    \left|\frac{1}{\zeta(\sigma+it)}\right| \le 1\,980 \frac{(\log{t})^{11/12}}{ (\log\log{t})^{3/4} }.
\end{equation}
Third, if 
\begin{equation}\label{eqn:KV-obj1}
    \sigma \ge 1-\frac{1}{53.99 (\log t)^{2/3}(\log\log t)^{1/3}} ,
    \quad\text{then}\quad
    \left|\frac{1}{\zeta(\sigma+it)}\right| \le 6\,368 (\log{t})^{2/3} (\log\log{t})^{1/4} .
\end{equation}
Subsequently, for $x \geq e^{4\,589.20}$, $x \geq e^{10\,352.03}$, and $x \geq \exp(e^{10.01})$, respectively, we have 
\begin{align}
    |M(x)| &< 33.56 x (\log{x}) e^{- \sqrt{\frac{\log{x}}{5.56}}} , \label{eqn:classical-obj2-hl} \\
    |M(x)| &< 21.48 x (\log{x}) e^{- \sqrt{\left(1 - \frac{\log{21.24}}{\log\log{x}}\right) \frac{\log{x} \log\log{x}}{42.48}}} , \quad\text{and} \label{eqn:Littlewood-obj2-hl} \\
    |M(x)| &< 25.85 x (\log{x}) e^{- \left(\frac{5}{3 (53.99)^3} \left(1 - \frac{\log{53.99}}{\log\log{x}}\right)^{-1}\right)^{1/5} \frac{(\log{x})^{3/5}}{(\log\log{x})^{1/5}}} . \label{eqn:KV-obj2-hl}
\end{align}
\end{theorem}

\begin{remark}
Prior to the recent result of the second author \cite[Cor.~5]{nicol}, the best bounds obtainable for the reciprocal $\zeta(\sigma+it)^{-1}$ within the zero-free regions \eqref{zfr regions} were of orders $\log{t}$, $\frac{\log{t}}{\log\log{t}}$, and $(\log{t})^{2/3}(\log\log{t})^{1/3}$, respectively. Therefore, the methods developed in the present work not only yield asymptotically sharper bounds for $\zeta(\sigma+it)^{-1}$ throughout these regions, but the explicit computations underlying Theorem \ref{thm:headline} also improve upon earlier results of \cite{CullyHugillLeong, nicol, Trudders} by orders of magnitude. 
For example, \cite[Cor.~4]{nicol} shows that if $t \ge 13$ and \begin{equation}\label{eqn:Leong} \sigma \ge 1-\frac{1}{5.56\log t}, \quad\text{then}\quad \left|\frac{1}{\zeta(\sigma+it)}\right| \le 7.59 \cdot 10^{31}\log t. \end{equation} Repeating the computations underlying \eqref{eqn:Leong} under the assumption $t\ge \mathcal H$ reduces the constant $7.59\cdot 10^{31}$ only slightly, to $7.13\cdot 10^{31}$. Even with this refinement, however, the resulting bound remains much weaker than \eqref{eqn:classical-obj1} of our Theorem~\ref{thm:headline}.
\end{remark}

\begin{remark}
Unconditional and explicit bounds of the form $M(x) \ll x/(\log{x})^k$ for integers $k \geq 0$ also exist in the literature (e.g., see \cite{chirre2025optimalboundssumsbounded, cde, sunconjarxiv2024, Ramare13}). One such result from Johnston, the second author, and Tudzi \cite[Thm.~A.1]{sunconjarxiv2024} states
\begin{equation}\label{eqn:JLT}
    | M(x) | \leq \frac{2.9189 x}{(\log{x})^2} \qquad \text{for all}\qquad x > 1.
\end{equation}
Bounds like these are particularly useful for small values of $x$. However, any explicit version of the bounds in \eqref{eqn:explicitise_me} would eventually be stronger as $x$ increases. Prior to the present work, explicit bounds for $M(x)$ of the forms in \eqref{eqn:explicitise_me} were not available in the published literature.
\end{remark}

\begin{remark}
It is expected that the computations underlying Theorem \ref{thm:headline} should hold on much broader ranges of $x$. In practice, though, the computed ranges in Theorem \ref{thm:headline} (and Theorem \ref{thm:objective2} in Section \ref{sec:results}) are more than sufficient. To illustrate this, one can compute thresholds $y_A < y_B < y_C$ such that \eqref{eqn:JLT} is (numerically) stronger than \eqref{eqn:classical-obj2-hl} for $\log{x} < y_A$, \eqref{eqn:classical-obj2-hl} is stronger than \eqref{eqn:Littlewood-obj2-hl} for $y_A \leq \log{x} < y_B$, \eqref{eqn:Littlewood-obj2-hl} is stronger than \eqref{eqn:KV-obj2-hl} for $y_B \leq \log{x} < y_C$, and \eqref{eqn:KV-obj2-hl} in turn is strongest for $\log{x} \geq y_C$. Rounding to two decimal places, we find that $y_A = 1\,650.77$, $y_B = 138\,529.33$, and $y_C = e^{22.54}$ are admissible, so by extension one can see \eqref{eqn:classical-obj2-hl}, \eqref{eqn:Littlewood-obj2-hl}, and \eqref{eqn:KV-obj2-hl} each hold for every $x \geq e^{4\,589.20}$. 
\end{remark}

\section{Results}\label{sec:results}

The ultimate goal of this paper is to prove unconditional and explicit versions of the bounds in \eqref{eqn:explicitise_me}. To this end, we first derive explicit upper bounds for $|\zeta(\sigma+it)|^{-1}$ within the zero-free regions described in \eqref{zfr regions} (Objective 1), which are likely to be of independent interest, and subsequently apply these bounds to obtain corresponding unconditional and explicit bounds for $M(x)$ (Objective 2). 

\subsection{Objective 1}

The first goal of this paper is to prove unconditional and explicit upper bounds for $\zeta(\sigma + it)^{-1}$ of the asymptotic orders given in \eqref{eqn:asymptotic_forms}. These bounds hold inside their respective forms of zero-free region given in \eqref{zfr regions}. Our results, presented below, significantly improve recent explicit bounds established in \cite{CullyHugillLeong, nicol, Trudders}, both numerically and asymptotically. 

\begin{theorem}\label{thm:objective1}
Let $t \geq \mathcal{H}$ with $\mathcal{H}$ defined as in \eqref{eqn:PlattTrudgian}. Let $\mathcal{R}_1$, $\mathcal{R}_2$, $\mathcal{R}_3$ be the corresponding values to fixed choices of $W$ in Table \ref{tab:objectives}. 
First, if 
\begin{equation}\label{eqn:classical-obj1-xtra}
    \sigma \ge 1-\frac{1}{W \log{t}} ,
    \quad\text{then}\quad
    \left|\frac{1}{\zeta(\sigma+it)}\right|\le \mathcal{R}_1 (\log{t})^{11/12} .
\end{equation}
Second, if 
\begin{equation}\label{eqn:Littlewood-obj1-xtra}
    \sigma \ge 1-\frac{\log\log{t}}{W \log{t}} ,
    \quad\text{then}\quad
    \left|\frac{1}{\zeta(\sigma+it)}\right|\le \mathcal{R}_2 \frac{(\log{t})^{11/12}}{ (\log\log{t})^{3/4}}.
\end{equation}
Third, if 
\begin{equation}\label{eqn:KV-obj1-xtra}
    \sigma \ge 1-\frac{1}{W (\log t)^{2/3}(\log\log t)^{1/3}} ,
    \quad\text{then}\quad
    \left|\frac{1}{\zeta(\sigma+it)}\right|\le \mathcal{R}_3 (\log{t})^{2/3} (\log\log{t})^{1/4} .
\end{equation}
\end{theorem}

\begin{table}[]
    \centering
    \resizebox{\textwidth}{!}{%
    \setlength{\tabcolsep}{3pt}
    \begin{tabular}{|cccc||cccc||cccc|}
    \hline
    $W$ & $\log{x_1(W)}$ & $\mathcal{R}_1$ & $c_1(x_1(W))$ & $W$ & $\log{x_2(W)}$ & $\mathcal{R}_2$ & $c_2(x_2(W))$ & $W$ & $\log\log{x_3(W)}$ & $\mathcal{R}_3$ & $c_3(x_3(W))$ \\
    \hline
    $5.56$ & $4\,589.20$ & $1\,237$ & $33.56$ & $21.24$ & $10\,352.03$ & $1\,980$ & $21.48$ & $53.99$ & $10.01$ & $6\,368$ & $25.85$ \\
    $6$ & $4\,952.38$ & $793$ & $22.21$ & $21.5$ & $10\,490.06$ & $942$ & $13.08$ & $54$ & $10.01$ & $5\,223$ & $22.19$ \\
    $7$ & $5\,777.77$ & $701$ & $18.16$ & $22$ & $10\,755.87$ & $790$ & $11.73$ & $55$ & $10.03$ & $488$ & $7.02$ \\
    $8$ & $6\,603.17$ & $686$ & $16.34$ & $25$ & $12\,360.29$ & $731$ & $10.57$ & $60$ & $10.11$ & $473$ & $6.85$ \\
    $10$ & $8\,253.96$ & $678$ & $14.06$ & $30$ & $15\,066.77$ & $721$ & $9.66$ & $70$ & $10.27$ & $467$ & $6.63$ \\
    \hline
    \end{tabular}%
    }
    \caption{The constants in Theorem \ref{thm:objective1} at fixed choices of $W$ and $i\in\{1,2,3\}$.}
    \label{tab:objectives}
\end{table}

\begin{remark}
For \eqref{eqn:classical-obj1-xtra} it is a requirement that $W > Z_1$, for \eqref{eqn:Littlewood-obj1-xtra} it is a requirement that $W > Z_2$, and for \eqref{eqn:KV-obj1-xtra} it is a requirement that $W > Z_3$, where $Z_1$, $Z_2$, and $Z_3$ are defined as in \eqref{z zerofree const}. 
In addition, the results \eqref{eqn:classical-obj1}, \eqref{eqn:Littlewood-obj1}, \eqref{eqn:KV-obj1} in our Theorem \ref{thm:headline} are a special case of Theorem \ref{thm:objective1} and the bounds \eqref{eqn:Littlewood-obj1} and \eqref{eqn:KV-obj1} in our Theorem \ref{thm:objective1} are the first explicit results of their kind. 
\end{remark}

\begin{remark} 
The constants $\mathcal{R}_1$, $\mathcal{R}_2$, and $\mathcal{R}_3$ are tied to both the choice of $W$ and the underlying zero-free region. The entries in Table \ref{tab:objectives} should therefore be read row-by-row within each block. For example, the bound \eqref{eqn:Littlewood-obj1-xtra} with $\mathcal{R}_2 = 1\,980$ corresponds to $W=21.24$, not $W=5.56$. 
\end{remark}

The proof of Theorem \ref{thm:objective1} combines explicit bounds for the logarithmic derivative $\tfrac{\zeta'(\sigma+it)}{\zeta(\sigma+it)}$ of the Riemann zeta function with a new computational smoothing procedure. In this context, computational smoothing refers to replacing a single large integration step across a zero-free region by many thousands (or even millions) of much smaller steps. More precisely, we perform a large collection of locally optimised computations involving explicit bounds for $\tfrac{\zeta'(\sigma+it)}{\zeta(\sigma+it)}$. By accumulating these intermediate computations rather than relying on one global bound, we substantially reduce the influence of the largest contributions, which arise near the boundary of the zero-free region. This significantly strengthens the resulting bounds. 

A limitation of our computational smoothing approach is that the computations required to establish Theorem \ref{thm:objective1} can be expensive. However, numerical strength may be traded for efficiency by reducing the number of intermediate calculations. In turn, this generally yields weaker bounds, although the resulting loss in numerical strength is mitigated in regions farther from the boundary of the underlying zero-free region.

Capitalising on this trade-off, we have proved the special cases \eqref{eqn:classical-obj1}, \eqref{eqn:Littlewood-obj1}, and \eqref{eqn:KV-obj1} of Theorem \ref{thm:objective1} using more extensive calculations, and then established the remaining bounds in Theorem \ref{thm:objective1} via a less computationally demanding approach.  Consequently, it still produces strong numerical results while substantially reducing the computational cost.

\subsection{Objective 2}

The second goal of this paper is to apply the bounds in Theorem \ref{thm:objective1} to obtain unconditional and explicit upper bounds for $|M(x)|$ of the asymptotic forms in \eqref{eqn:explicitise_me}. Our results are given in Theorem \ref{thm:objective2} below. Note that $\log{x} \ll x^{\varepsilon}$ for any $\varepsilon > 0$, and so this factor can be absorbed into $\eta_1, \eta_2, \eta_3$ in \eqref{eqn:explicitise_me} if desired. 

\begin{theorem}\label{thm:objective2}
For fixed choices of $W$ and $i\in\{1,2,3\}$, let $\mathcal{R}_i$ and $x_i(W)$ be the corresponding values in Table \ref{tab:objectives}. Also let
\begin{align}
    c_1(x) 
    &= 2e + \frac{12 \mathcal{R}_1}{23\pi W} \left(\frac{\log{x}}{W}\right)^{-1/24}
    + \frac{4 e^{2 + \frac{1}{\log{x}}}}{\sqrt{W\log{x}}} , \label{eqn:c1_def} \\
    c_2(x) 
    &= 2e + \left( 1 + \frac{3}{4\log\log{\mathcal{H}}}\right) \frac{12 \mathcal{R}_2  (\log\log{\mathcal{H}})^{1/4}}{23\pi W (\log{\mathcal{H}})^{1/12}}
    + 4 e^{2+\frac{1}{\log{x}}} \sqrt{\frac{\log\log{x}}{W\log{x}}} , \quad\text{and}\label{eqn:c2_def} \\
    c_3(x) 
    &= 2e + \frac{3 \mathcal{R}_3}{5 \pi W(\log\log{\mathcal{H}})^{1/12}}
    + \frac{4 e^{2+\frac{1}{\log{x}}}}{W^{3/5} (\log{x})^{2/5}} . \label{eqn:c3_def}
\end{align} 
We have $c_i(x) \leq c_i(x_i(W))$ where values for $c_i(x_i(W))$ are recorded in Table \ref{tab:objectives}, 
\begin{align}
    |M(x)| &< c_1(x,W,\mathcal{R}_1) x (\log{x}) e^{- \sqrt{\frac{\log{x}}{W}}} \quad\text{for}\quad x \geq x_1(W) , \label{eqn:classical-ob2} \\
    |M(x)| &< c_2(x,W,\mathcal{R}_2) x (\log{x}) e^{- \sqrt{\left(1 - \frac{\log{W}}{\log\log{x}}\right) \frac{\log{x} \log\log{x}}{2 W}}} \quad\text{for}\quad x \geq x_2(W), \quad\text{and} \label{eqn:Littlewood-ob2} \\
    |M(x)| &< c_3(x,W,\mathcal{R}_3) x (\log{x}) e^{- \left(\frac{5}{3 W^3} \left(1 - \frac{\log{W}}{\log\log{x}}\right)^{-1}\right)^{1/5} \frac{(\log{x})^{3/5}}{(\log\log{x})^{1/5}}} \quad\text{for}\quad x \geq x_3(W) . \label{eqn:KV-ob2} \\
\end{align}
\end{theorem}

\begin{remark}
Meagre refinements are possible for the results \eqref{eqn:Littlewood-ob2} and \eqref{eqn:KV-ob2}. That is, by re-tracing our arguments, one finds that $c_2(x) = c_3(x) = 2e + o(1)$ is admissible as $x\to\infty$ in Theorem \ref{thm:objective2}, provided one uses the appropriate lower bound from Lemma \ref{lem:logT_bounds} in place of the cruder bound $T \geq \mathcal{H}$. 
We did not pursue this refinement, as it would complicate the definitions of $c_2$, $c_3$, while yielding only a marginal improvement.
\end{remark}

\begin{remark}
The bound \eqref{eqn:classical-ob2} in Theorem \ref{thm:objective2} is a significant refinement of earlier work by Chalker \cite{chalker}. 
For example, if $\log{x} \geq 79\,589.39$ and $W = 6$, then Chalker (see \cite[Tab.~3.3]{chalker}) proved $c_1(x) < 1.7\cdot 10^{29}$ is admissible in \eqref{eqn:classical-ob2}. For the same choice of $W$, our computations in Theorem \ref{thm:objective2} certify $c_1(x) \leq 22.21$ in a larger range. This level of improvement is both feasible and justified, because our approach refines the standard contour-shifting argument and utilises more efficient computational techniques.
\end{remark}

\subsection{Computations}\label{ssec:computations}

The computations in this paper have been executed using the \texttt{Python} programming language on a standard personal computer. Our code is freely available on GitHub \cite{LeeLeongCode}. This repository contains code sufficient to reproduce the computations underlying Table \ref{tab:objectives} and Theorems \ref{thm:objective1}-\ref{thm:objective2}. 

\subsection{Structure of the paper}

The remainder of this paper is devoted to the proofs of our main results, namely Theorems \ref{thm:objective1}-\ref{thm:objective2}. In Section \ref{sec:bounds_LDRZF}, we develop a unified framework for obtaining explicit bounds for the logarithmic derivative of $\zeta(s)$ in any of the zero-free regions described in \eqref{zfr regions}. Section \ref{sec:bounds_RRZF} uses this framework to establish Theorem \ref{thm:objective1}. Finally, in Section \ref{sec:bounds_MF}, we combine these results with additional arguments to prove Theorem \ref{thm:objective2}.

To streamline the proofs presented in this paper, we introduce several parameters. While these may initially appear somewhat technical, we hope to preserve readability by ensuring that all parameters introduced in a given section are local to that section and are used only within its context.

\section{Parametrised bounds for the logarithmic derivative}\label{sec:bounds_LDRZF}

To prepare for the arguments that follow and the statement of Proposition \ref{prop:parametrised_bounds_ld}, we first introduce some notation. Unless indicated otherwise, we adhere to this notation throughout this section. Suppose there exists $\xi(t)$ such that for real $2/3 \leq \xi_1 \leq 1$ and $-1 \leq \xi_2 \leq 1/3$, we have 
\begin{equation}\label{eqn:xi_def}
    (\log{t})^{\xi(t)} = (\log{t})^{\xi_1} (\log\log{t})^{\xi_2} \leq \log{t} .
\end{equation}
Of particular interest are the choices
\begin{equation*}
    (\xi_1,\xi_2) = (1,0) , \quad
    (\xi_1,\xi_2) = (1,-1) , \quad\text{and}\quad
    (\xi_1,\xi_2) = (\tfrac{2}{3},\tfrac{1}{3}) ,
\end{equation*}
as they recover the denominator in the classical, Littlewood, and Korobov--Vinogradov zero-free regions respectively. 

For any $w > 0$, we also define the decreasing functions
\begin{equation}\label{a2(t0)}
    a_{w}(t) : = 1+
    \frac{\log \left ( 1+\frac{w}{t}\right )}{\log t}
    \quad\text{and}\quad 
    b_{w}(t) := 1+ \frac{\log{a_{w}(t)}}{\log \log t} .
\end{equation}
In our argument, we will construct concentric discs with origin $s_{0} = \sigma_{0} +it$, where
\begin{equation}\label{sigma0}
    \sigma_{0} = 1 + \frac{d}{(\log{t})^{\xi_1} (\log\log{t})^{\xi_2}} .
\end{equation}
Further, for some $\omega\ge 0$, we let the radius of the outermost disc be
\begin{equation}\label{r def}
    \begin{split}
        r &= \frac{d + \omega \log\log{t}}{(\log{t})^{\xi_1} (\log\log{t})^{\xi_2}} .
    \end{split}
    \end{equation}

In this section, we extend and update the methods of \cite{CullyHugillLeong, hiaryleongyangArxiv, nicol} to prove a unified framework from which bounds can be established for the logarithmic derivative of $\zeta(\sigma + it)$ in the classical, Littlewood, and Korobov--Vinogradov zero-free regions. In particular, we prove Proposition \ref{prop:parametrised_bounds_ld}, from which one can compute explicit bounds for the logarithmic derivative of $\zeta(s)$ in each of the zero-free regions described in \eqref{zfr regions}. 

\begin{proposition}\label{prop:parametrised_bounds_ld}
Fix any $W > 0$ such that $\zeta (\sigma +it)\neq 0$ in the region
\begin{equation}\label{eqn:zfr_condition}
    \sigma \ge 1-\frac{1}{W (\log{t})^{\xi_1} (\log\log{t})^{\xi_2}}
    \quad\text{and}\quad 
    t\ge t_{0}\ge \mathcal{H} .
\end{equation}
In addition, we let: 
\begin{itemize}
    \item $d\leq Z^{-1}$ such that $Z$ is the smallest value $W$ can take; 
    \item $\omega \geq \frac{W^{-1}}{\log\log{t}}$ be a constant such that $B_\omega(t) \ll 1$ and
    \begin{equation}\label{eqn:omega1}
         \omega \leq \frac{(\log{t})^{\xi_1}}{2 (\log\log t)^{1 - \xi_2}} ,
    \end{equation} 
    where $B_{\omega}(t)$ is defined in \eqref{zeta_AB-B} below;
    \item $r$ be defined as in \eqref{r def}; 
    \item $a_r(t)$, $b_r(t)$ be defined as in \eqref{a2(t0)}; and
    \begin{equation}
        \beta = \left(\frac{d^{-1}}{W} + 1\right) \left( a_{r}(t)^{-\xi_1} b_{r}(t)^{-\xi_2} + 1 \right)^{-1} \label{beta}.
        \end{equation}
    
\end{itemize}
Suppose also that the chosen parameters satisfy the condition
\begin{equation}\label{eqn:alpha_cond}
    \frac{(\omega / d) \log\log{t} + 1}{2(a_{r}(t)^{-\xi_1} b_{r}(t)^{-\xi_2} + 1)} >1 .
\end{equation} 
Then in the region \eqref{eqn:zfr_condition}, we have
\begin{equation*}
    \left| \frac{\zeta'(\sigma +it)}{\zeta(\sigma +it)} \right| \le Q_{d, \omega, W, t_0} (\log t)^{\xi_1}(\log\log t)^{\xi_2},
\end{equation*}
in which
\begin{align*}
    Q_{d, \omega, W, t_0} 
    &= \max \Bigg \lbrace 
    \Omega(\beta, d, \omega, t_0), \left(\frac{1}{W} + 2d\right)^{-1} 
    \Bigg\rbrace
\end{align*}
with $A_4(d,\omega,t,r)$ defined as in \eqref{A4} and
\begin{align*}
    \Omega(\beta, d, \omega, t)  
    &= \left( \frac{8\beta}{1-\beta} \right) 
    \left( B_\omega(t) + \frac{3\xi_1}{4} + \frac{1}{6} + \frac{3\xi_2 \log\log\log{t}}{4\log\log{t}} + \frac{\log{A_4(d,\omega,t,r)}}{\log\log{t}} \right) \\
    &\qquad\qquad\qquad\cdot\left(1 + \frac{d}{\omega\log\log{t}}\right) \left(1 - d \frac{1 + 2 a_{r}(t)^{-\xi_1} b_{r}(t)^{-\xi_2}}{\omega\log\log{t}}\right)^{-1} \\
    &\qquad\qquad\qquad\cdot\left(\omega - d \frac{1 + 2 a_{r}(t)^{-\xi_1} b_{r}(t)^{-\xi_2}}{\log\log{t}}\right)^{-1} \\
    &\qquad\qquad + \left( \frac{1+\beta}{1-\beta} \right) \frac{1}{d} .
\end{align*}
Furthermore, if $\sigma \geq 1$, we can replace \eqref{beta} with $\beta = \left( a_{r}(t)^{-\xi_1} b_{r}(t)^{-\xi_2} + 1 \right)^{-1}$, and then replace the definition of $Q_{d, \omega, W, t_0}$ with 
\begin{align*}
    Q_{d, \omega, W, t_0} &= \max \Bigg \lbrace 
    \Omega(\beta, d, \omega, t_0), 
    \frac{1}{2d} 
    \Bigg\rbrace .
\end{align*}
\end{proposition}

The remainder of this section is devoted to the proof of Proposition \ref{prop:parametrised_bounds_ld}. We begin by collecting some auxiliary technical results that are central to the argument in Section \ref{ssec:ab}. In Section \ref{ssec:tech_ingredients}, we establish two further technical lemmas (Lemmas \ref{lem:A1}-\ref{lem:A2}) that allow us to apply these results with the required parameters asserted. Finally, in Section \ref{ssec:proof_prop}, we combine these ingredients to complete the proof. 

\subsection{Auxiliary results}\label{ssec:ab}

For convenience, we collect some important auxiliary results here.

\begin{lemma}[\text{\cite[Lem.~15]{hiaryleongyangArxiv}}]\label{log_deriv_zeta_lem1}
Let $s_{0}$ be a complex number. Let $r$ and $\alpha $ be positive numbers (possibly depending on $s_{0}$) such that $\alpha <1/2$. Suppose that the function $f(s)$ is analytic in a region containing the disc $|s-s_{0}|\le r$. Suppose further there is a number $A_{1}$ independent of $s$ such that
\begin{equation*}
    \left |\frac{f(s)}{f(s_{0})}\right | \le A_{1}
    \quad\text{for}\quad
    |s-s_{0}|\le r.
\end{equation*}
Then, for any $s$ in the disc $|s-s_{0}|\le \alpha r$ we have
\begin{equation*}
\left | \frac{f'(s)}{f(s)} - \sum _{\rho }\frac{1}{s-\rho}\right |
\le \frac{4\log A_{1}}{r(1 -2\alpha )^{2}},
\end{equation*}
where $\rho $ runs through the zeros of $f(s)$ in the disc
$|s -s_{0}|\le \tfrac{r}{2}$, counted with multiplicity.
\end{lemma}

\begin{lemma}[\text{\cite[Lem.~16]{hiaryleongyangArxiv}}]\label{lem_3.2_HB}
Let $s$ and $s_{0}$ be complex numbers with real parts $\sigma $ and
$\sigma _{0}$, respectively. Let $r$, $\alpha $, $\beta $, $A_{1}$, and
$A_{2}$ be positive numbers, possibly depending on $s_{0}$, such that
$\alpha <1/2$ and $\beta <1$. Suppose that the function $f(s)$ satisfies
the conditions of Lemma~\ref{log_deriv_zeta_lem1} with $r$,
$\alpha$, and $A_{1}$, and that
\begin{equation*}
\left | \frac{f'(s_{0})}{f(s_{0})}\right | \le A_{2}.
\end{equation*}
Suppose, in addition, that $f(s) \neq 0$ for any $s$ in the intersection
of the disc $|s-s_{0}|\le r$ and the right half-plane
$\sigma \ge \sigma _{0} - \alpha r$. Then, for any $s$ in the disc
$|s-s_{0}|\le \alpha \beta r$,
\begin{equation*}
\left | \frac{f'(s)}{f(s)}\right | \le
\frac{8\beta \log A_{1}}{r(1-\beta )(1-2\alpha )^{2}} +
\frac{1+\beta}{1-\beta}A_{2}.
\end{equation*}
\end{lemma}

We also prove the following results, which are important for our later arguments. 

\begin{lemma}\label{zeta_log_bound}
Let $\omega \ge 0$ be given. Suppose $t \ge t_0$,
\begin{equation}\label{omega assumption}
    \sigma_t := 1 - \frac{\omega\log\log t}{(\log{t})^{\xi_1} (\log\log{t})^{\xi_2}} ,
    \quad \text{and} \quad 
    \omega \leq \frac{(\log{t})^{\xi_1}}{2 (\log\log t)^{1 - \xi_2}}.
\end{equation}
Then, if $\sigma \ge \sigma_t$, we have
\begin{equation*}
    |\zeta(\sigma + it)| \le A({\omega},t_0)(\log t)^{B_\omega(t)} ,
\end{equation*}
in which
\begin{align}
    A({\omega},t_0) 
    &:=
    \begin{cases}
        70.6995 &\text{if $\omega >0$, $t_0 \ge e^e$,}\\
        58.096\sqrt{1 + 9t_0^{-2}}\left(1 + \frac{\log \sqrt{1 + 9t_0^{-2}}}{\log t_0}\right) &\text{if $\omega =0$, $t_0 \ge 3$,}
    \end{cases}
    \quad\text{and}\label{zeta_AB-A} \\
    B_\omega(t) 
    &:=
    \frac{2}{3} + 4.43795 \omega^{3/2} (\log\log{t})^{\frac{1 - 3\xi_2}{2}} (\log t)^{1 - \frac{3\xi_1}{2}} . \label{zeta_AB-B}
\end{align}
\end{lemma}

\begin{proof}
The proof is similar to that of \cite[Lem.~7]{hiaryleongyangArxiv}, and is split into three cases.

\medskip\noindent\textsc{Case I:} 
To begin, we prove the result for $\sigma < 1$. 
It follows from \cite[Thm.~1.1]{bellotti2024explicit} that if $1/2 \le \sigma \le 1$ and $t\ge 3$, then
\begin{equation}\label{ford bound}
    |\zeta(\sigma + it)| \le 70.6995 t^{4.43795(1-\sigma)^{3/2}}(\log t)^{2/3}.
\end{equation}
Since $\sigma \geq \sigma_t \ge 1/2$ is asserted, we require
\begin{equation*}
    \omega \leq \frac{(\log{t})^{\xi_1}}{2 (\log\log t)^{1 - \xi_2}}
\end{equation*}
to ensure \eqref{ford bound} can be applied. 
Note also that $(1-\sigma)^{3/2}$ decreases monotonically as $\sigma$ increases. So, if $\sigma_t \le \sigma \le 1$ and $\omega > 0$, then \eqref{ford bound} gives
\begin{align*}%\label{expr_exp}
    |\zeta(\sigma + it)| 
    &\le 70.6995\exp\left(4.43795(1-\sigma_t)^{3/2}\log t + \frac{2}{3}\log\log t\right) \\
    &\le 70.6995\exp\left(4.43795\left( \frac{\omega\log\log t}{(\log{t})^{\xi_1}(\log\log t)^{\xi_2}} \right)^{3/2}\log t + \frac{2}{3}\log\log t\right) \\
    &= A({\omega},t_0) (\log{t})^{B_\omega(t)} ,
\end{align*}
where $A({\omega},t_0)$ and $B_\omega(t)$ are defined as in \eqref{zeta_AB-A}-\eqref{zeta_AB-B}. With this, the result is proved for every $\sigma_t \le \sigma < 1$, because $\omega = 0$ implies $\sigma_t = 1$.

\medskip\noindent\textsc{Case II:} 
Next, we prove the result holds in the range $1 \le \sigma \le 2$ using the Phragm\'en--Lindel\"of Principle \cite[Lem. 3]{nicol} with the holomorphic function $f(s) = (s - 1)\zeta(s)$.  
Assuming this choice, we bound $f(\sigma + it)$ on the line $\sigma = 1$ and $\sigma = 2$ separately. To this end, recall from \cite[Lem. 5]{nicol} that for $|t| \ge 3$, we have
\begin{equation*}
    |\zeta(1+it)| \le 58.096(\log t)^{2/3} .
\end{equation*}
It follows that
\begin{equation*}
    |f(1 + it)| = |t|\,| \zeta(1 + it)| \le 58.096\, |2 + it|(\log{|2 + it|})^{2/3} 
\end{equation*}
and
\begin{equation*}
    |f(2 + it)| \le |1 + it|\, \zeta(2) < 58.096|3 + it|(\log |3 + it|)^{2/3} .
\end{equation*}
The first inequality is verified numerically for $|t| < 3$, so these inequalities hold for every real $t$. 
So, \cite[Lem. 3]{nicol} with $a = 1$, $b = 2$, $Q = 1$, $\alpha_1 = \beta_1 = 1$, $\alpha_2=\beta_2 = 2/3$, and $A = B = 58.096$ tells us that
\begin{equation}\label{pl bound}
    |f(s)| \le 58.096|1 + s|(\log |1 + s|)^{2/3}
\end{equation}
for $1 \le \sigma \le 2$. For every $t \geq t_0$, we also have
\begin{equation*}
    \begin{split}
        \frac{|1 + s|}{|s - 1|} &\le \frac{\sqrt{9 + t^2}}{t} \le \sqrt{1 + 9t_0^{-2}}
        \quad\text{and}\\
        \log|1 + s| &\le \log\sqrt{9 + t^2} \le \log t + \log\sqrt{1 + 9t_0^{-2}} .
    \end{split}
\end{equation*}
Thus for $1 \le \sigma \le 2$, the inequality \eqref{pl bound} implies the desired result; this argument holds for any choice of $\omega\geq 0$.

\medskip\noindent\textsc{Case III:} 
Finally, we address the case where $\sigma \ge 2$. In this case, we have
\begin{equation*}
    |\zeta(\sigma+it)|\le \zeta(2) \le 1.65 \le A(0,t_0) (\log t)^{2/3} \le A({\omega},t_0) (\log t)^{B_\omega(t)}
\end{equation*}
for all $t\ge t_0$. 
\end{proof}

\begin{proposition}\label{lb 1/z s>1}
Let $d_1>0$ be a real parameter and $\delta_1 := d_1 (\log t)^{-\xi_1} (\log\log t)^{-\xi_2}$. Then for $\sigma \ge 1+\delta_1$ and $t\ge t_0 \ge 3$, we have
\begin{equation}
    \left|\frac{1}{\zeta(\sigma+it)}\right| 
    \le A_3(d_1,t_0) (\log{t})^{\frac{3\xi_1}{4} + \frac{1}{6}}(\log\log t)^{\frac{3\xi_2}{4}} ,
\end{equation}
with
\begin{equation}\label{A3}
    A_3(d_1,t_0) := \frac{A(0,t_0)^{1/6}}{d_1^{3/4}} \left( 1+\frac{\log 2}{\log t_0}\right)^{1/6} \exp\left(\frac{3\gamma_0 d_1/4}{(\log t_0)^{2/3}(\log\log t_0)^{1/3}}\right),
\end{equation}
where $\gamma_0$ is Euler's constant and $A({\omega},t_0)$ is defined in Lemma \ref{zeta_log_bound}.
\end{proposition}

\begin{proof}
Recall $2(1+\cos{\theta})^2 = 3 + 4\cos{\theta} + \cos{2\theta} \geq 0$ for any real $\theta$. Following the classical non-negativity argument with this observation (see \cite[\S 3.3]{titchmarsh1986theory} and \cite{leongmossinghoffArxiv} for more commentary on this approach), we have $\zeta^3(\sigma)|\zeta^4(\sigma+it)\zeta(\sigma +2it)| \ge 1$ for $\sigma > 1$, which subsequently implies
\begin{equation}\label{trig ineq}
    \left| \frac{1}{\zeta(1+\delta_1+it)}\right| \le |\zeta(1+\delta_1)|^{3/4}|\zeta(1+\delta_1 + 2it)|^{1/4} .
\end{equation}
By Lemma \ref{zeta_log_bound} with $\omega = 0$, we have
\begin{align}
    |\zeta(1+\delta_1 +2it)|^{1/4} \le A(0,t_0)^{1/6} (\log (2t))^{1/6} \le A(0,t_0)^{1/6}\left( 1+\frac{\log 2}{\log t_0}\right)^{1/6}\left( \log t\right)^{1/6} . \label{trig bound 2}
\end{align}
Further, for any $\sigma >1$, \cite[Lem.~5.4]{ramare2016explicit} states that $\zeta(\sigma) \le e^{\gamma_0 (\sigma - 1)}/\sigma - 1$, so
\begin{equation}
    |\zeta(1+\delta_1)|^{3/4} 
    \le \left(\frac{(\log t)^{\xi_1}(\log\log t)^{\xi_2}}{d_1}\right)^{3/4} \exp\bigg(\frac{3\gamma_0 d_1/4}{(\log t)^{\xi_1}(\log\log t)^{\xi_2}}\bigg) .
\end{equation}
Inserting these upper bounds in \eqref{trig ineq} proves the result.
\end{proof}

\subsection{Technical ingredients}\label{ssec:tech_ingredients}

Fix $\xi_1$, $\xi_2$ as in \eqref{eqn:xi_def} and $t\ge t_{0} \ge \mathcal{H}$. Let $r$ and $\sigma_0$ be defined as in \eqref{sigma0}-\eqref{r def}. 
In order to apply Lemmas \ref{log_deriv_zeta_lem1}-\ref{lem_3.2_HB} later, we require the following intermediate results (Lemmas \ref{lem:A1}-\ref{lem:A2}). Note that Lemma \ref{lem:A1} is a consequence of Lemma \ref{zeta_log_bound} and Proposition \ref{lb 1/z s>1}. 

\begin{lemma}\label{lem:A1}
Suppose $\omega$ satisfies \eqref{eqn:omega1} and recall the definitions of $r$, $A$, $B_{\omega}$, $A_3$ from \eqref{r def}, \eqref{zeta_AB-A}, \eqref{zeta_AB-B}, \eqref{A3} respectively. Let $s_0 = \sigma_0 +it$. For every $s$ in the disc $|s - s_0| \leq r$, we have  
\begin{align*}
|\Im{s} - t| \leq r,\qquad
    \Re{s} &\ge 1 - \frac{\omega (\log\log{t})^{1-\xi_2}}{(\log{t})^{\xi_1} } , 
\end{align*}
and for $t\ge t_0$, we have
\begin{align*}
    \left |\frac{\zeta (s)}{\zeta (\sigma_{0}+it)} \right| 
    &\le A_{1} := A_4(d,\omega,t_0,r) (\log{t})^{B_\omega(t) + \frac{3\xi_1}{4} + \frac{1}{6}}(\log\log t)^{\xi_2 + \frac{3\xi_2}{4}} ,
\end{align*}
    where
    \begin{equation}\label{A4}
        A_4(d,\omega,t,r) = \left(1+\frac{r}{t}\right)^{4.43795 \left(\frac{\omega (\log\log{t})^{1-\xi_2}}{(\log{t})^{\xi_1}}\right)^{3/2}} a_r(t)^{B_\omega(t)} 
        A_3(d,t) A(\omega,t-r) .
    \end{equation}
\end{lemma}

\begin{proof}
Observe that points in the disc $|s - s_0| \leq r$ satisfy
\begin{equation*}
    s 
    = (\sigma_0 + r_1\cos\theta) + i(t + r_1\sin\theta)
    := u(\theta) + i v(\theta) ,
\end{equation*}
where $0\leq r_1 \leq r$ and $\theta \in (0,2\pi]$ is the argument. 
For every $s$ satisfying
\begin{equation*}
    u(\theta) \geq 1 - \frac{\omega (\log\log(t+r))^{1-\xi_2}}{(\log(t+r))^{\xi_1}}
    \quad\text{or}\quad
    v(\theta) \leq t ,
\end{equation*}
it follows from Lemma \ref{zeta_log_bound} that
\begin{align*}
    |\zeta(s)| 
    &\leq A(\omega,t_0-r) (a_r(t) \log{t})^{B_\omega(t+r)}
    \leq A(\omega,t_0-r) (a_r(t) \log{t})^{B_\omega(t)} .
\end{align*} 
Upon combination with Proposition \ref{lb 1/z s>1}, this yields a result that is marginally stronger than the bound in the statement of the lemma. 

All that remains is to prove the result for $s$ which satisfy $\theta \in [\pi/2,\pi)$,
\begin{equation*}
    1 - \frac{\omega (\log\log{t})^{1-\xi_2}}{(\log{t})^{\xi_1}} 
    \leq u(\theta)
    \leq 1 - \frac{\omega (\log\log(t+r))^{1-\xi_2}}{(\log(t+r))^{\xi_1}} ,
\end{equation*}
and $t < v(\theta) \leq t + r$. 
It follows from the condition \eqref{eqn:omega1} that $u(\theta) \geq 1/2$ in this case. Therefore, we can apply \eqref{ford bound} and the above bounds to see
\begin{align*}
    |\zeta(s)| 
    &\leq 70.6995 (t+r)^{4.43795(1- u(\theta))^{3/2}}(a_r(t)\log{t})^{2/3} \\ 
    &\leq 70.6995 \left(t+r\right)^{4.43795 \left(\frac{\omega (\log\log{t})^{1-\xi_2}}{(\log{t})^{\xi_1}}\right)^{3/2}} (a_r(t) \log{t})^{2/3} \\ 
    &= \left(1+\frac{r}{t}\right)^{4.43795 \left(\frac{\omega (\log\log{t})^{1-\xi_2}}{(\log{t})^{\xi_1}}\right)^{3/2}} a_r(t)^{2/3} 
    A(\omega,t_0) (\log{t})^{B_\omega(t)} .
\end{align*} 
Combine the above with Proposition \ref{lb 1/z s>1} to yield the result.
\end{proof}

\begin{lemma}\label{lem:A2}
For $s_0$ defined as in \eqref{sigma0}, we have
\begin{equation}\label{A2}
    \left| \frac{\zeta'(s_0)}{\zeta(s_0)}\right| < A_2 : =\frac{(\log{t})^{\xi_1} (\log\log{t})^{\xi_2}}{d}.
\end{equation}
\end{lemma}

\begin{proof}
For any $\sigma >1$, by \cite{delange1987remarque} we have
    \begin{equation}\label{A2 trivial}
    \left| \frac{\zeta'(\sigma +it)}{\zeta(\sigma+it)}\right| \le -\frac{\zeta'(\sigma)}{\zeta(\sigma)} < \frac{1}{\sigma -1} .
    \end{equation}
Applying this bound to our choice of $s_0$ yields the result.
\end{proof}

\subsection{Proof of Proposition \ref{prop:parametrised_bounds_ld}}\label{ssec:proof_prop}
We aim to apply Lemma \ref{lem_3.2_HB}. Suppose $\zeta(s) \neq 0$ for every $s$ in the intersection $I$ of the disc $|s-s_{0}|\le r$ and the right half-plane $\sigma \ge \sigma_{0} - \alpha r$ for some $\alpha < 1/2$. 
Note that $\sigma \geq \sigma_0 - \alpha r$ implies
\begin{align*} 
    \sigma 
    \geq \sigma_0 - \frac{\alpha (d+\omega\log\log t)}{(\log{t})^{\xi_1} (\log\log{t})^{\xi_2}} 
    &= 1 + \frac{(1 - \alpha) d - \alpha\omega\log\log t}{(\log{t})^{\xi_1} (\log\log{t})^{\xi_2}} \\
    &= 1 - \left( \alpha\omega + \frac{(\alpha - 1) d}{\log\log{t}} \right) \left(\frac{(\log\log{t})^{1-\xi_2}}{(\log{t})^{\xi_1} }\right) .
\end{align*}
Therefore, to ensure the entirety of $I$ is inside the zero-free region described in \eqref{eqn:zfr_condition}, it suffices to assert $d \leq Z^{-1}$ and 
\begin{equation*}
    \left( \alpha\omega + \frac{(\alpha - 1) d}{\log\log{t}} \right) \left(\frac{(\log\log{t})^{1-\xi_2}}{(\log{t})^{\xi_1} }\right) 
    \leq \frac{d}{(\log t_r)^{\xi_1} (\log\log t_r)^{\xi_2}}, 
\end{equation*}
in which $t_r = t + r$. It follows from $\log{t_r} = a_r(t)\log{t}$, $\log\log{t_r} = b_r(t)\log\log{t}$, and the above observations that we require
\begin{align}
    \alpha 
    &\leq \frac{d \left( a_{r}(t)^{-\xi_1} b_{r}(t)^{-\xi_2} + 1 \right)}{\omega \log\log{t} + d} . \label{eqn:alpha_ub}
\end{align}
We also require $\alpha < 1/2$, which follows from the assumption \eqref{eqn:alpha_cond}. 
Now, we can apply Lemma \ref{lem_3.2_HB} with $f(s) = \zeta(s)$ and $A_{1}$, $A_2$ defined as in Lemmas \ref{lem:A1}-\ref{lem:A2}. For $\alpha$ satisfying the condition \eqref{eqn:alpha_ub} and suitable $\beta < 1$, in the disc $|s-s_{0}|\le \alpha \beta r$, it follows that
\begin{align*}
    \left| \frac{\zeta'(s)}{\zeta(s)} \right| 
    &\le \frac{8\beta \log{A_{1}}}{r(1-\beta )(1-2\alpha )^{2}} 
    + \frac{1+\beta}{1-\beta}A_{2} \\
    &= \left( \frac{8\beta}{r(1-\beta )(1-2\alpha )^{2}} \right) \\
    &\qquad\qquad\cdot\left( B_\omega(t) + \frac{3\xi_1}{4} + \frac{1}{6} + \frac{3\xi_2 \log\log\log{t}}{4\log\log{t}} + \frac{\log{A_4(d,\omega,t_0,r)}}{\log\log{t}} \right) \log\log{t} \\
    &\qquad\quad + \left( \frac{1+\beta}{1-\beta} \right) \frac{(\log{t})^{\xi_1} (\log\log{t})^{\xi_2}}{d} .
\end{align*}
From the definition \eqref{r def} and the inequality \eqref{eqn:alpha_ub}, we also have 
\begin{align*}
    r (1-2\alpha)^{2}
    &\geq \frac{d + \omega \log\log{t}}{(\log{t})^{\xi_1} (\log\log{t})^{\xi_2}} \left(\frac{\omega\log\log{t} - d (1 + 2 a_{r}(t)^{-\xi_1} b_{r}(t)^{-\xi_2})}{\omega\log\log{t} + d}\right)^2 \\
    &= \frac{(\log\log{t})^{1-\xi_2}}{(\log{t})^{\xi_1}} \left(\frac{1 - d \frac{1 + 2 a_{r}(t)^{-\xi_1} b_{r}(t)^{-\xi_2}}{\omega\log\log{t}}}{1 + \frac{d}{\omega\log\log{t}}}\right) 
    \left(\omega - d \frac{1 + 2 a_{r}(t)^{-\xi_1} b_{r}(t)^{-\xi_2}}{\log\log{t}}\right) .
\end{align*}
Therefore, the upper bound becomes
\begin{align*}
    \left| \frac{\zeta'(s)}{\zeta(s)} \right| 
    &\le \Omega(\beta, d, \omega, t_0) (\log{t})^{\xi_1} (\log\log{t})^{\xi_2} .
\end{align*}
In practice, it is optimal to have $\alpha$ as large as possible, so we let $\alpha$ equal the upper bound in \eqref{eqn:alpha_ub}. Therefore, we have now proved bounds on $|\zeta'(s)/\zeta(s)|$ which hold for any $s = \sigma + it$ such that
\begin{equation}\label{eqn:sig_holds_for}
\begin{split}
    1 - \frac{d \left( \beta \left( a_{r}(t)^{-\xi_1} b_{r}(t)^{-\xi_2} + 1 \right) - 1\right)}{(\log t)^{\xi_1}(\log\log t)^{\xi_2}}
    &\leq \sigma 
    \leq 1 + \frac{d \left( \beta \left( a_{r}(t)^{-\xi_1} b_{r}(t)^{-\xi_2} + 1 \right) + 1\right)}{(\log t)^{\xi_1}(\log\log t)^{\xi_2}} .
\end{split}
\end{equation}
To ensure the lower bound in \eqref{eqn:sig_holds_for} matches the boundary of the region \eqref{eqn:zfr_condition}, we select
\begin{equation}\label{eqn:W3_condition}
    \beta = \left(\frac{d^{-1}}{W} + 1\right) \left( a_{r}(t)^{-\xi_1} b_{r}(t)^{-\xi_2} + 1 \right)^{-1} .
\end{equation}
To prove an upper bound on $|\zeta'(s)/\zeta(s)|$ that holds for $\sigma$ larger than the RHS of \eqref{eqn:sig_holds_for}, we apply \eqref{A2 trivial} and \eqref{eqn:W3_condition} to reveal
\begin{align*}
    \left| \frac{\zeta'(s)}{\zeta(s)}\right| 
    < \frac{1}{\sigma -1}
    &\leq \frac{(\log t)^{\xi_1}(\log\log t)^{\xi_2}}{d \left( \beta \left( a_{r}(t)^{-\xi_1} b_{r}(t)^{-\xi_2} + 1 \right) + 1\right)} \\
    &= \left(\frac{1}{W} + 2d\right)^{-1} (\log t)^{\xi_1}(\log\log t)^{\xi_2} .
\end{align*}
Finally, taking the maximum of these upper bounds yields the result. 

The case for $\sigma \ge 1$ follows exactly the same argument, using the different choice of $\beta$ stated. \qed

\section{Bounds for the reciprocal of the Riemann zeta function}\label{sec:bounds_RRZF}

Using Proposition \ref{prop:parametrised_bounds_ld}, we can compute explicit upper bounds for the logarithmic derivative in each of the zero-free regions described in \eqref{zfr regions}. From these computations, we can subsequently prove strong bounds for the reciprocal of $\zeta(s)$ in each of the zero-free regions described in \eqref{zfr regions} using the upcoming result (Proposition \ref{thm:1/zeta 3 inter}). While our proof of Proposition \ref{thm:1/zeta 3 inter} is quite straightforward, we must first establish one preparatory bound, namely Lemma \ref{lem:prep_for_inter}, which follows from Proposition \ref{prop:parametrised_bounds_ld}.

\begin{lemma}\label{lem:prep_for_inter}
Let $\sigma \geq 1$, $t \geq t_0 \geq \mathcal{H}$, and
\begin{equation*}
    (\xi_1,\xi_2,\xi_C) \in \left\{(1,0, 24.303) ,\,(1,-1,110.6),\,\left(\tfrac{2}{3},\tfrac{1}{3},216.667\right) \right\} ,
\end{equation*}
then we have 
\begin{align*}
    \left|\frac{\zeta'(\sigma+it)}{\zeta(\sigma+it)}\right| \le \xi_C (\log t)^{\xi_1}(\log\log t)^{\xi_2}.
\end{align*}
\end{lemma}

\begin{proof}
The first case, where $(\xi_1,\xi_2) = (1,0)$ is proved in \cite[Cor. 2]{nicol}, and the second case, where $(\xi_1,\xi_2) = (1,-1)$, is proved in \cite[Tab. 4]{CullyHugillLeong}. All that remains is to prove the remaining third case.

In the case $\sigma \geq 1$, Proposition \ref{prop:parametrised_bounds_ld} provides a stronger definition of $Q_{d, \omega, W, t_0}$ for all $t \geq t_0 \geq \mathcal{H}$, and \cite[Cor. 2]{nicol} also tells us that for $\sigma \ge 1$ and $13 \le t \le t_0$, we have
\begin{align*}
    \left|\frac{\zeta'(\sigma+it)}{\zeta(\sigma+it)}\right| 
    %&\le 24.303\left( \frac{\log t}{\log\log t}\right)^{1/3}(\log t)^{2/3}(\log\log t)^{1/3} \\
    &\le 24.303\left( \frac{\log t_0}{\log\log t_0}\right)^{1/3} (\log{t})^{2/3} (\log\log{t})^{1/3} . 
\end{align*}
So, we set $\log{t_0} = 5000$ and increment this choice until the largest $\log{t_0}$ (rounded to two decimal places) is located such that $24.303 (\log t_0 / \log\log t_0 )^{1/3} < Q_{d, \omega, W, t_0}$. 
In the end, we find that $\log{t_0} = 6186.07$ is the largest such $t_0$, where 
\begin{equation*}
     \max\left\{ Q_{d, \omega, W, t_0}, \, 24.303\left( \frac{\log t_0}{\log\log t_0}\right)^{1/3}  \right\} < 216.667 . \qedhere
\end{equation*}
\end{proof}

\begin{proposition}\label{thm:1/zeta 3 inter}
Let there be real parameters $d_1 >0$, $t_{0} \ge \mathcal{H}$, and $(\xi_1,\xi_2,\xi_C)$ as in Lemma \ref{lem:prep_for_inter}.
Let $\{ V_{j} : 1\leq j\leq J\}$ be a sequence of increasing real numbers where $(V_{j}, Q_{H,j})$ is a pair for which $t \ge t_{0}$ and
\begin{equation*}
    \sigma \ge 1- \frac{1}{V_{j}(\log t)^{\xi_1}(\log\log t)^{\xi_2}} \quad \text{implies}\quad \left |
    \frac{\zeta'(\sigma +it)}{\zeta(\sigma +it)}\right| \le Q_{H,j}(\log t)^{\xi_1}(\log\log t)^{\xi_2}.
\end{equation*}
Then, if $t\ge t_0$ and
\begin{equation*}
    \sigma \ge 1- \frac{1}{V_{1}(\log t)^{\xi_1}(\log\log t)^{\xi_2}},
    \quad\text{we have}\quad
    \left | \frac{1}{\zeta (\sigma +it)}\right | \le R (\log t)^{\frac{3\xi_1}{4} + \frac{1}{6}}(\log\log t)^{\frac{3\xi_2}{4}} ,
\end{equation*}
where
\begin{align*}
    R &= A_3(d_1,t_0)\cdot\max \Biggl\{ 1,\,  \exp \Bigg( \sum _{j=1}^{J-1} Q_{H,j} \left( \frac{1}{V_{j}}-\frac{1}{V_{j+1}}\right) +
    \frac{Q_{H,J}}{V_{J}} +\xi_C d_{1} \Bigg) \Biggr\rbrace , 
\end{align*}
with $A_3(d_1,t_0)$ defined as in \eqref{A3}. 
In addition, if $\sigma \geq 1$, then $R$ can be replaced with $$R = A_3(d_1,t_0)\cdot \max\{ 1, \exp(\xi_C d_{1})\}.$$
\end{proposition}

\begin{proof}
Let $\delta _{1} = \delta _{1}(t) := d_{1}(\log t)^{-\xi_1}(\log\log t)^{-\xi_2}$. If $\sigma > 1 + \delta_1$, then the result follows from Proposition \ref{lb 1/z s>1}. So, all that remains is to prove the result for $\sigma \le 1+\delta _{1}$, and then take the maximum of the two. 
To this end, note that in the range $1- (V_1(\log t)^{\xi_1}(\log\log t)^{\xi_2})^{-1}\le \sigma \le 1+ \delta _{1}$ we have 
\begin{align}\label{1/zeta_int}
\log \left| \frac{1}{\zeta (\sigma +it)}\right| = - \Re
\log \zeta \left( 1+ \delta_{1} +it\right) + \int_{\sigma}^{1+\delta_{1}} \Re\left( \frac{\zeta'}{\zeta}(y+it)\right) \,dy.
\end{align}
Writing $\Delta = (\log t)^{-\xi_1}(\log\log t)^{-\xi_2}$, we can split the integral and rewrite it as
\begin{align}\label{smooth}
    \left ( \sum_{j = 1}^{J-1} \int _{1-\tfrac{\Delta}{V_{j}}}^{1-\tfrac{\Delta}{V_{j+1}}} + \int _{1-\tfrac{\Delta}{V_{J}}}^{1} + \int _{1}^{1+\delta _{1}}
    \right ) \Re\left ( \frac{\zeta '}{\zeta}(y+it)\right )
    \,dy.
\end{align}
It follows from Lemma \ref{lem:prep_for_inter} and our assumptions that
\begin{align*}
    \int_{\sigma}^{1+\delta _{1}} \Re\left ( \frac{\zeta '}{\zeta}(y+it)\right ) \,dy
    &\leq \sum _{j=1}^{J-1} Q_{H,j}
\left ( \frac{1}{V_{j}}-\frac{1}{V_{j+1}}\right ) +
\frac{Q_{H,J}}{V_{J}} +\xi_C d_{1}.
\end{align*}
Note that if $\sigma \geq 1$, then in the above upper bound we can assert 
\begin{align*}
    \sum_{j=1}^{J-1} Q_{H,j}
    \left( \frac{1}{V_{j}}-\frac{1}{V_{j+1}} \right) +
    \frac{Q_{H,J}}{V_{J}} = 0.
\end{align*}
Therefore, from \eqref{1/zeta_int} we have
\begin{equation}\label{1/zeta_int1}
\begin{split}
    \log \left| \frac{1}{\zeta (\sigma +it)}\right| 
    &< -\log \left |
    \zeta \left ( 1+\delta_1+ it \right )\right| \\
    &\qquad + \sum _{j=1}^{J-1} Q_{H,j} \left ( \frac{1}{V_{j}} - \frac{1}{V_{j+1}}\right ) + \frac{Q_{H,J}}{V_{J}} +\xi_C d_{1} .
\end{split}
\end{equation}
To bound the first term on the right-hand side of \eqref{1/zeta_int1}, apply Proposition~\ref{lb 1/z s>1}. This gives
\begin{equation}\label{V_term_VTEX1}
    \left | \frac{1}{\zeta (1 + \delta_{1} +it)}\right | \le A_3(d_1,t_0) (\log t)^{\frac{3\xi_1}{4} + \frac{1}{6}}(\log\log t)^{\frac{3\xi_2}{4}}.
\end{equation}
Combine this with \eqref{1/zeta_int1} and exponentiate to yield the result for $\sigma \le 1+\delta _{1}$. %With this, the proof is complete. 
\end{proof}

Next, we prove Theorem \ref{thm:objective1} by making suitable choices in Proposition \ref{thm:1/zeta 3 inter} as follows. 

\begin{proof}[Proof of Theorem \ref{thm:objective1}]
To begin, we prove the result in the Korobov--Vinogradov zero-free region. For each fixed choice of $W$ in Table \ref{tab:objectives}, we set $t_0 = \mathcal{H}$ and integers $k$ such that 
\begin{equation}\label{subdivisions-xtra}
    V(W) 
    %= \{ V_{j} : 1\leq j\leq J\}
    = \{ W + k\cdot 10^{-5} : k\in [0,2\cdot 10^{5})\} \cup \{ W + 0.5 k : k\in [4,1000] \} .
\end{equation}
In the special case $W = 53.99$, we replace \eqref{subdivisions-xtra} with
\begin{equation}\label{subdivisions}
    V(W) 
    = \{ W + k\cdot 10^{-6} : k\in [0,2\cdot 10^{6})\} \cup \{ W + 0.5 k : k\in [4,1000] \} .
\end{equation}
Next, for each $v\in V(W)$, we compute $Q_v$ using Proposition \ref{prop:parametrised_bounds_ld} such that $t\geq t_0$ and
\begin{equation*}
    \sigma \ge 1- \frac{1}{v (\log t)^{2/3}(\log\log t)^{1/3}} 
    \quad \text{imply}\quad \left |
    \frac{\zeta'(\sigma +it)}{\zeta (\sigma +it)}\right | \le Q_{v}(\log t)^{2/3}(\log\log t)^{1/3} .
\end{equation*}
To compute each of the values $Q_v$, we assert $\xi_1 = 2/3$, $\xi_2 = 1/3$, $\xi_C = 216.667$, $Z = Z_3$ defined as in \eqref{z zerofree const}, $t_0 = \mathcal{H}$ defined as in \eqref{eqn:PlattTrudgian}, $d = Z^{-1}$, and $\omega = 1.5046945258$ in Proposition \ref{prop:parametrised_bounds_ld}. This choice of $\omega$ is almost-optimal for each $W = v$ under consideration (up to ten decimal places); we have determined this by systematically incrementing $\omega$ until the upper bound in Proposition \ref{prop:parametrised_bounds_ld} stops improving on a broad selection of $W$. Carrying out these calculations requires an amount of computation time that is determined by the volume of $v \in V(W)$, but the resulting improvement in our values for $R$ more than justify the effort. 

Next, we set $d_1 = 10^{-8}$ and increment this choice appropriately until the definition of $R$ (as given in Proposition \ref{thm:1/zeta 3 inter}) is minimised; invariably the result is optimised (to eight decimal places) upon asserting $d_1 = 0.00346104$. Therefore, in light of the above computations for $Q_v$ and this choice, it is straightforward to compute the values $\mathcal{R}_3 = R$ associated to each choice of $W$ in Table \ref{tab:objectives}.  

To prove the result in the classical zero-free region, we follow an analogous process. That is, for each fixed choice of $W$ in Table \ref{tab:objectives}, we assert $\xi_1 = 1$, $\xi_2 = 0$, $\xi_C = 24.303$, $Z = Z_1$ defined as in \eqref{z zerofree const}, $t_0 = \mathcal{H}$, $d = Z^{-1}$, $\omega = 3.4812610455$ in Proposition \ref{prop:parametrised_bounds_ld} to compute values $Q_v$ for every $v\in V(W)$ such that 
$t\geq t_0$ and
\begin{equation*}
    \sigma \ge 1- \frac{1}{v \log{t}} 
    \quad \text{imply}\quad \left |
    \frac{\zeta'(\sigma +it)}{\zeta (\sigma +it)}\right | \le Q_{v}\log{t} .
\end{equation*}
In this case, we define $V(W)$ as in \eqref{subdivisions-xtra} for $W > 5.56$, while for $W = 5.56$, we use
\begin{equation*}
    V(W) 
    = \{ W + k\cdot 10^{-6} : k\in [0,2\cdot 10^{7})\} \cup \{ W + 0.5 k : k\in [4,1000] \} .
\end{equation*}
Apply Proposition \ref{thm:1/zeta 3 inter} with $d_1 = 0.00049$ and these computations to recover the values for $\mathcal{R}_1$ associated to each choice for $W$ presented in Table \ref{tab:objectives}. 

Finally, to prove the result in the Littlewood zero-free region for each fixed choice of $W$ in Table \ref{tab:objectives}, we assert $\xi_1 = 1$, $\xi_2 = -1$, $\xi_C = 110.6$, $Z = Z_2$ defined as in \eqref{z zerofree const}, $t_0 = \mathcal{H}$, $d = Z^{-1}$, $\omega = 1.0144550824$ in Proposition \ref{prop:parametrised_bounds_ld} to compute values $Q_v$ for every $v\in V(W)$ such that 
$t\geq t_0$ and
\begin{equation*}
    \sigma \ge 1- \frac{\log\log{t}}{v \log{t}} 
    \quad \text{imply}\quad \left |
    \frac{\zeta'(\sigma +it)}{\zeta (\sigma +it)}\right | \le Q_{v}\frac{\log{t}}{\log\log{t}} .
\end{equation*}
In this case, we define $V(W)$ as in \eqref{subdivisions-xtra} for $W > 21.24$, and $V(W)$ as in \eqref{subdivisions} for $W = 21.24$. 
Apply Proposition \ref{thm:1/zeta 3 inter} with $d_1 = 0.00069$ and these computations to recover the values for $\mathcal{R}_2$ associated to each choice for $W$ presented in Table \ref{tab:objectives}. 
\end{proof}

In addition to Theorem \ref{thm:objective1}, to bound $M(x)$ we also require bounds for the $\zeta(\sigma+it)^{-1}$ that hold in the region $|t| \leq \mathcal{H}$ of the critical strip. In this region, it is known that $\zeta(\sigma + it) \neq 0$ in the half plane $\sigma > 1/2$, so results over broader regions are available. 
In the slightly restricted range $\sigma \geq 3/4$, we prove the following result (Lemma \ref{lem:central_inequality_lwr}), which suffices for our purposes. 

We also recall from Simoni\v{c} \cite[Lem.~4]{Simonic22} that in the restricted range $1/2 < \sigma \leq 3/2$ and $|t| \leq 2\exp(e^2)$, we have
\begin{equation}\label{eqn:SimonicLem4}
    \frac{1}{|\zeta(\sigma + it)|} \leq \frac{4}{\sigma - \frac{1}{2}} .
\end{equation}

\begin{lemma}\label{lem:central_inequality_lwr}
For $\sigma \geq 0.75$ and $14 \leq t\leq \mathcal{H}$, we have
\begin{align*}
    \left| \frac{1}{\zeta(\sigma + it)} \right| 
    \leq 29.388t^{4.964} \log t 
    &\leq e^{146} \log{t} ,
\end{align*}
where the implied constant is maximised at $\mathcal{H}$. 
\end{lemma}

\begin{proof}
First recall from \cite[Cor. 4]{nicol} that for any $\sigma \geq 1$ and $t \geq 13$ we have
\begin{align*}%\label{eqn:central_inequality_uno}
    \left| \frac{1}{\zeta(\sigma + it)} \right| 
    \leq 30.812\log t
    < 29.388t^{4.964} \log t .
\end{align*} 
Hence all that remains is to prove the result for $1 > \sigma \ge 0.75$ and $14 \le t \leq \mathcal{H}$. To this end, we have
\begin{equation}\label{eqn:central_inequality_dos}
    \left| \frac{1}{\zeta(\sigma + it)} \right| 
    = \exp\left( - \Re{\log{\zeta(1 + it)}} + \int_{\sigma}^{1} \Re{\frac{\zeta'(y+it)}{\zeta(y+it)}}\,dy \right) .
\end{equation}
Now, recall from \cite[Cor. 4]{nicol} that for $t\ge 2$, we have the bound
\begin{align*}
    \exp\big(- \Re{\log{\zeta(1 + it)}} \big)
    = \left| \frac{1}{\zeta(1 + it)} \right| \le 29.388\log t.
\end{align*}
In addition, \cite[Cor.~3,~Tab.~2]{nicol} implies that
\begin{align*}
    \int_{\sigma}^{1} \Re{\frac{\zeta'(y+it)}{\zeta(y+it)}}\,dy 
    &\leq \sum_{(\sigma_1,\sigma_2)\in S_1} (\sigma_2 - \sigma_1) \max_{\sigma_1 
    \leq y \leq \sigma_2} \left| \frac{\zeta'(y+it)}{\zeta(y+it)} \right| \\
    &\leq \left(\frac{35.162 + 23.759 + 17.050 + 13.123 + 10.175}{20}\right) \log{t}\\
    &< 4.964 \log{t} 
\end{align*}
when summing over the pairs $$(\sigma_1,\sigma_2) \in S_1 = \{(0.75,0.8), (0.8,0.85), (0.85,0.9), (0.9,0.95), (0.95,1)\}.$$ 
Inserting these upper bounds into \eqref{eqn:central_inequality_dos} completes the proof.
\end{proof}

\section{Bounds for Mertens function}\label{sec:bounds_MF}

In this section, we prove Theorem \ref{thm:objective2} using Theorem \ref{thm:objective1}. 
To begin, we establish a common framework in Section \ref{ssec:setup} that underpins each proof, before proving each bound \eqref{eqn:classical-ob2}, \eqref{eqn:Littlewood-ob2}, and \eqref{eqn:KV-ob2} in Theorem \ref{thm:objective2} in Sections \ref{ssec:classical}, \ref{ssec:Littlewood}, and \ref{ssec:KV} respectively. 

\subsection{Common framework}\label{ssec:setup}

Suppose $F(s) = \sum_{n\geq 1} a_n n^{-s}$ denotes a Dirichlet series, $c > 0$ is a real parameter larger than the abscissa of absolute convergence of $F$, $x> 1$, and $T\geq 1$. The truncated Perron formula tells us
\begin{equation}\label{eqn:ogPF}
    \left| \sum_{n \leq x} a_n - \frac{1}{2\pi i} \int_{c - iT}^{c + iT} F(s)\frac{x^s}{s}\,ds \right|
    \leq E(x,T,c) 
\end{equation}
for some error $E(x,T,c)$ dependent on $x$, $T$, and $c$. Simoni\v{c} applies \cite[Thm.~7.1]{RamareELSI} in \cite[Thm.~4]{Simonic22} to demonstrate that
\begin{equation}\label{eqn:RamarePF}
    E(x,T,c)
    = \frac{2 F(c) x^c}{T} + 4 e^c\left(\frac{e\phi(x) x\log{T}}{T} + \phi(ex)\right) 
\end{equation}
is admissible, where $\phi(n)$ satisfies $|a_n| \leq \phi(n)$. Stronger formulations of the error in \eqref{eqn:ogPF} are available in the literature (see e.g., \cite{CullyHugillJohnston, CullyHugillJohnston2}). However, we prefer to work with the simpler formulation \eqref{eqn:RamarePF}, since the improvements afforded by the stronger versions would be small relative to the additional technical considerations they would introduce.

In the special case $a_n = \mu(n)$, we have $\phi(n) = 1$ and $F(c) \leq \zeta(c) \leq c/(c-1)$, which converges absolutely on $\Re{s} > 1$. 
Therefore, it follows from \eqref{eqn:ogPF} and \eqref{eqn:RamarePF} that for any $x > 1$, $c > 1$, and some $T\geq 1$, we have
\begin{equation}\label{eqn:approximation}
    \left| M(x) - \frac{1}{2\pi i} \int_{c - iT}^{c + iT} \frac{x^s}{s \zeta(s)}\,ds \right|
    \leq \frac{2 c x^c}{T (c-1)} + 4 e^c\left(\frac{e x\log{T}}{T} + 1\right) .
\end{equation}

In Section \ref{ssec:pb_i}, we bound the contour integral in \eqref{eqn:approximation} by moving the line of integration inside the critical strip; our result is presented in Lemma \ref{lem:important_integral_Chalker}. In Section \ref{ssec:optimal_order}, we identify a balanced choice for the truncation point $T$. In Section \ref{ssec:completion}, we combine these observations to prove Proposition \ref{prop:MF_framework}, which establishes parametrised upper bounds for $|M(x)|$ from which we prove our main results. 

\subsubsection{Parametrised bounds for the contour integral in \eqref{eqn:approximation}}\label{ssec:pb_i}

Note the notation $f_1$, $f_2$, in the following result is designed to ensure the framework is general enough to support our proof of \eqref{eqn:classical-ob2}, \eqref{eqn:Littlewood-ob2}, and \eqref{eqn:KV-ob2}.  

\begin{lemma}\label{lem:important_integral_Chalker}
Suppose $0.75 \leq \sigma_0 < 1$ and $T \geq \mathcal{H}$, where $\mathcal{H}$ is defined as in \eqref{eqn:PlattTrudgian}. Also, let
\begin{equation*}
    (\lambda_1,\lambda_2,\lambda_3,\lambda_4) \in \left\{ (1,0,\tfrac{11}{12},0), (1,-1,\tfrac{11}{12},-\tfrac{3}{4}), (\tfrac{2}{3},\tfrac{1}{3},\tfrac{2}{3},\tfrac{1}{4})\right\}
\end{equation*}
such that 
\begin{equation*}
    f_1(t) = (\log{t})^{\lambda_1} (\log\log{t})^{\lambda_2}
    \qquad\text{and}\qquad 
    f_2(t) = (\log{t})^{\lambda_3} (\log\log{t})^{\lambda_4} .
\end{equation*} 
Further, let $W > 0$ be such that for any $t \geq \mathcal{H}$, in the region 
\begin{equation*}
    \sigma \ge \sigma_W(t) := 1-\frac{1}{W f_1(t)} 
\end{equation*}
we have $\zeta(\sigma+it)\neq 0$ and $\left|\zeta(\sigma+it)\right|^{-1} \le R_1 f_2(t)$ for some constant $R_1 > 0$. For any $\sigma \geq \sigma_0$ and $14 \leq t \leq \mathcal{H}$, we also suppose that 
\begin{equation*}
    \left|\frac{1}{\zeta(\sigma+it)}\right|\le R_2 \log{t} ,
\end{equation*}
for some constant $R_2 > 0$. Then, for $x> 1$, $c = 1 + 1/\log{x}$, and $\sigma_0 \leq \sigma_W(\mathcal{H})$, we have
\begin{equation*}
    \left|\frac{1}{2\pi i} \int_{c-iT}^{c+iT} \frac{x^s}{s\zeta(s)}\,ds \right| 
    \leq \mathcal{E}_{R_1,R_2,\sigma_0}(x,T) ,
\end{equation*}
where $\delta_{\lambda_4} = 1$ if $\lambda_4 < 0$, $\delta_{\lambda_4} = 0$ if $\lambda_4 \geq 0$, and
\begin{align*}
    \mathcal{E}_{R_1,R_2,\sigma_0}(x,T)
    &= \frac{R_1}{\pi \log{x}} \left( \frac{e x (\log{T})^{\lambda_3} (\log\log{T})^{\lambda_4}}{T} + \frac{x^{\sigma_W(\mathcal{H})} (\log{\mathcal{H}})^{\lambda_3} (\log\log{\mathcal{H}})^{\lambda_4}}{\mathcal{H}} \right) \\
    &\hspace{2cm} + \frac{R_1 x^{\sigma_W(T)}}{\pi} \Bigg( \left( 1 - \frac{\delta_{\lambda_4} \lambda_4}{\log\log{T}}\right) \frac{(\log{T})^{\lambda_3 + 1} (\log\log{T})^{\lambda_4}}{\lambda_3 + 1} \\ 
    &\hspace{8.6cm}- \frac{(\log{\mathcal{H}})^{\lambda_3 + 1} (\log\log{\mathcal{H}})^{\lambda_4}}{\lambda_3 + 1} \Bigg)  \\
    &\hspace{2cm} + \frac{x^{\sigma_0}}{\pi} \left( \frac{R_2 (\log{\mathcal{H}})^2}{2} + \frac{56}{\sigma_0(\sigma_0 - \tfrac{1}{2})} \right) .
\end{align*}
\end{lemma}

\begin{figure}
\centering
\begin{tikzpicture}[xscale=1.2,yscale=1.0]

\coordinate (A) at (7,-3.5);      % c-iT
\coordinate (B) at (7,3.5);       % c+iT
\coordinate (C) at (4.2,3.5);     % sigma_W(T)+iT
\coordinate (D) at (3.0,1.2);   % sigma_W(H)+iH
\coordinate (E) at (0,1.2);     % sigma_0+iH
\coordinate (F) at (0,-1.2);    % sigma_0-iH
\coordinate (G) at (3.0,-1.2);  % sigma_W(H)-iH
\coordinate (H) at (4.2,-3.5);    % sigma_W(T)-iT

\begin{scope}[very thick,
decoration={
 markings,
 mark=at position 0.5 with {\arrow{>}}
}]

% C0
\draw[postaction={decorate}]
(A)--(B)
node[midway,right=3pt] {$C_0$};

% C1
\draw[postaction={decorate}]
(B)--(C)
node[midway,above=3pt] {$C_1$};

% C2
\draw[postaction={decorate}]
(C) .. controls (3.7,2.4) and (4.0,1.8) ..
(D)
node[pos=0.55,right=3pt] {$C_2$};

% C3
\draw[postaction={decorate}]
(D)--(E)
node[midway,above=3pt] {$C_3$};

% C4
\draw[postaction={decorate}]
(E)--(F)
node[midway,left=3pt] {$C_4$};

% C5
\draw[postaction={decorate}]
(F)--(G)
node[midway,below=3pt] {$C_5$};

% C6
\draw[postaction={decorate}]
(G) .. controls (4.0,-1.8) and (3.7,-2.4) ..
(H)
node[pos=0.45,right=3pt] {$C_6$};

% C7
\draw[postaction={decorate}]
(H)--(A)
node[midway,below=3pt] {$C_7$};

\end{scope}

% Labels
\node[right] at (B) {$c+iT$};
\node[right] at (A) {$c-iT$};
\node[left] at (C) {$\sigma_W(T)+iT$};
\node[left] at (H) {$\sigma_W(T)-iT$};
\node[below] at (D) {$\sigma_W(\mathcal H)+i\mathcal H$};
\node[above] at (G) {$\sigma_W(\mathcal H)-i\mathcal H$};
\node[above] at (E) {$\sigma_0+i\mathcal H$};
\node[below] at (F) {$\sigma_0-i\mathcal H$};

\end{tikzpicture}
\caption{A rough illustration of the closed contour $C$ in the proof of Lemma \ref{lem:important_integral_Chalker}.}
\label{fig:contour}
\end{figure}

\begin{proof}
Define the contours $C_j$ such that
\begin{itemize}
    \item $C_0$ connects the nodes $c-iT$ and $c+iT$ via a straight line,
    \item $C_1$ connects the nodes $c+iT$ and $\sigma_W(T)+iT$ via a straight line,
    \item $C_2$ connects the nodes $\sigma_W(T)+iT$ and $\sigma_W(\mathcal{H})+i\mathcal{H}$ along the path $\sigma_W(t) + it$,
    \item $C_3$ connects the nodes $\sigma_W(\mathcal{H})+i\mathcal{H}$ and $\sigma_0 + i\mathcal{H}$ via a straight line,
    \item $C_4$ connects the nodes $\sigma_0+i\mathcal{H}$ and $\sigma_0 - i\mathcal{H}$ via a straight line,
    \item $C_5$ connects the nodes $\sigma_0-i\mathcal{H}$ and $\sigma_W(\mathcal{H})-i\mathcal{H}$ via a straight line,
    \item $C_6$ connects the nodes $\sigma_W(\mathcal{H})-i\mathcal{H}$ and $\sigma_W(T)-iT$ along the path $\sigma_W(t) - it$,
    \item $C_7$ connects the nodes $\sigma_W(T)-iT$ and $c - iT$ via a straight line.
\end{itemize}
These contours are purposefully chosen to avoid the non-trivial zeros $\varrho = \beta + i\gamma$ of $\zeta(s)$, through the use of the Riemann height and zero-free regions. In fact, the closed contour $$C := \bigcup_{j=0}^{7} C_j$$ contains no zeros of $\zeta(s)$, and hence Cauchy's residue theorem tells us
\begin{equation*}
    I 
    = \frac{1}{2\pi i} \int_{C} \frac{x^s}{s\zeta(s)}\,ds 
    = 0 .
\end{equation*}
For convenience, we give a rough illustration of the contour $C$ in Figure \ref{fig:contour}, and write
\begin{equation*}
    I_j := \frac{1}{2\pi i} \int_{C_j} \frac{x^s}{s\zeta(s)}\,ds ,
\end{equation*}
for $0\le j \le 7$, from which it follows that
\begin{equation*}
    \frac{1}{2\pi i} \int_{c-iT}^{c+iT} \frac{x^s}{s\zeta(s)}\,ds = - \sum_{j=1}^{7} I_j .
\end{equation*}
By our assumptions, we have
\begin{align*}
    %\sum_{j\in\{1,7\}} |I_j|
    %&\leq \frac{R_1}{\pi} \frac{f_2(T)}{T} \int_{\sigma_W(T)}^{c} x^{y}\,dy
    |I_1|+|I_7|
    &\leq \frac{R_1}{\pi} \frac{(\log{T})^{\lambda_3} (\log\log{T})^{\lambda_4}}{T} \int_{\sigma_W(T)}^{c} x^{y}\,dy 
    \leq \frac{R_1 e}{\pi} \frac{x (\log{T})^{\lambda_3} (\log\log{T})^{\lambda_4}}{T \log{x}} ,
\end{align*}
and
\begin{align*}
    |I_3|+|I_5|
    &\leq \frac{R_1}{\pi} \frac{(\log{\mathcal{H}})^{\lambda_3}(\log\log{\mathcal{H}})^{\lambda_4}}{\mathcal{H}} \int^{\sigma_W(\mathcal{H})}_{\sigma_0} x^{y}\,dy 
    \leq \frac{R_1}{\pi} \frac{(\log{\mathcal{H}})^{\lambda_3}(\log\log{\mathcal{H}})^{\lambda_4}}{\mathcal{H}} \frac{x^{\sigma_W(\mathcal{H})}}{\log{x}} .
\end{align*}
Through integration by parts, we also have
\begin{align*}
    |I_2|+|I_6|
    &\leq \frac{R_1 x^{\sigma_W(T)}}{\pi} \int_{\mathcal{H}}^{T} \frac{(\log{t})^{\lambda_3} (\log\log{t})^{\lambda_4}}{t}\,dt \\
    &= \frac{R_1 x^{\sigma_W(T)}}{\pi} \left( \left[ \frac{(\log{t})^{\lambda_3 + 1} (\log\log{t})^{\lambda_4}}{\lambda_3 + 1}\right]_{\mathcal{H}}^{T} - \frac{\lambda_4}{\lambda_3 + 1} \int_{\mathcal{H}}^{T} \frac{(\log{t})^{\lambda_3} (\log\log{t})^{\lambda_4 - 1}}{t}\,dt \right) .
\end{align*}
If $\lambda_4 \geq 0$, then we can discard the latter term for a more practical upper bound. 
The only case where $\lambda_4 < 0$ satisfies $(\lambda_1,\lambda_2,\lambda_3,\lambda_4) = (1,-1,\tfrac{11}{12},-\tfrac{3}{4})$, and we can treat the extra positive contribution in this special case directly. To this end, we observe that $(\log{t})^{\lambda_3} (\log\log{t})^{\lambda_4 - 1}$ increases on $\mathcal{H} \leq t \leq T$, and hence
\begin{align*}
    \frac{1}{(\log{T})^{\lambda_3 + 1} (\log\log{T})^{\lambda_4}} \int_{\mathcal{H}}^{T} \frac{(\log{t})^{\lambda_3} (\log\log{t})^{\lambda_4 - 1}}{t}\,dt 
    &\leq \frac{1}{\log{T} \log\log{T}} \int_{\mathcal{H}}^{T} \frac{dt}{t} 
    \leq \frac{1}{\log\log{T}} .
\end{align*}
Therefore, 
\begin{align*}
    |I_2|+|I_6|
    &\leq \frac{R_1 x^{\sigma_W(T)}}{\pi} \Bigg( \left( 1 - \frac{\delta_{\lambda_4} \lambda_4}{\log\log{T}}\right) \frac{(\log{T})^{\lambda_3 + 1} (\log\log{T})^{\lambda_4}}{\lambda_3 + 1} \\
    &\hspace{8cm}- \frac{(\log{\mathcal{H}})^{\lambda_3 + 1} (\log\log{\mathcal{H}})^{\lambda_4}}{\lambda_3 + 1} \Bigg) .
\end{align*}
Finally, \eqref{eqn:SimonicLem4} and our assumptions imply
\begin{align*}
    |I_4|
    &\leq \frac{x^{\sigma_0}}{\pi} \left( \int_{14}^{\mathcal{H}} \frac{|\zeta(\sigma_0 + it)|^{-1}}{t}\,dt + \frac{1}{\sigma_0} \int_{0}^{14} |\zeta(\sigma_0 + it)|^{-1}\,dt \right) \\
    &\leq \frac{x^{\sigma_0}}{\pi} \left( R_2 \int_{14}^{\mathcal{H}} \frac{\log{t}}{t}\,dt  + \frac{4\cdot 14}{\sigma_0 (\sigma_0 - \tfrac{1}{2})} \right)
    \leq \frac{x^{\sigma_0}}{\pi} \left( \frac{R_2 (\log{\mathcal{H}})^2}{2}  + \frac{56}{\sigma_0 (\sigma_0 - \tfrac{1}{2})} \right) .
\end{align*}
Combine the triangle inequality and these observations to deduce the result. 
\end{proof}

\subsubsection{The point of truncation}\label{ssec:optimal_order}

For fixed choices of $W$, $\lambda_1$, and $\lambda_2$, regardless of $\lambda_3$ and $\lambda_4$, we choose $T = T_{x}$ satisfying $T_{x} \geq \mathcal{H}$ and
\begin{equation}\label{eqn:T_choice}
    (\log{T_x})^{\lambda_1 + 1} (\log\log{T_x})^{\lambda_2} = \frac{\log{x}}{W} .
\end{equation}
Importantly, the relationship \eqref{eqn:T_choice} certifies $x^{\sigma_W(T)} = x/T$, which simplifies our upcoming arguments. One could attempt to optimise the choice of $T$, but doing so would require balancing several competing terms and would introduce a disproportionate amount of technical and notational complexity. Moreover, we believe the choice \eqref{eqn:T_choice} is near optimal as it strikes a good balance among the terms involved in the upper bound for $M(x)$.

\subsubsection{Completing the argument}\label{ssec:completion}

Asserting the choice for $T$ as in \eqref{eqn:T_choice}, we prove the following result, which provides parametrised bounds from which we can prove Theorem \ref{thm:objective2}.

\begin{proposition}\label{prop:MF_framework}
Assume the notation in Lemma \ref{lem:important_integral_Chalker} and fix admissible choices for $W$, $\lambda_1$, $\lambda_2$, $\lambda_3$, and $\lambda_4$. In addition, define $T = T_x$ satisfying \eqref{eqn:T_choice} and suppose that
\begin{align}\label{eqn:suffices_to_show}
\begin{split}
    &\frac{R_1 (\log{\mathcal{H}})^{\lambda_3 + 1} (\log\log{\mathcal{H}})^{\lambda_4}}{\pi (\lambda_3 + 1)} \\
    &\qquad\qquad> 2e + \frac{R_1}{\pi} \left( \frac{e (\log{\mathcal{H}})^{\lambda_3 - \lambda_1 - 1} (\log\log{\mathcal{H}})^{\lambda_4 - \lambda_2}}{W} + \frac{(\log{\mathcal{H}})^{\lambda_3} (\log\log{\mathcal{H}})^{\lambda_4}}{\mathcal{H} \log{x}} \right) \\
    &\qquad\qquad\qquad \,\, +\frac{\exp\left({\log{T} + (\sigma_0 - 1) \log{x}}\right)}{\pi} \left( \frac{R_2 (\log{\mathcal{H}})^2}{2} + \frac{56}{\sigma_0(\sigma_0 - \tfrac{1}{2})} + \frac{4 e^c \pi}{x^{\sigma_0}} \right)
\end{split}
\end{align}
for every $x \geq x_0$, where $x_0 > 0$ is the least value such that $T_x \geq \mathcal{H}$. 
Then for $x \geq x_0$, we have 
\begin{equation}\label{eqn:MF_framework}
\begin{split}
    \frac{| M(x) |}{x (\log{x}) /T} 
    &< 2e + \left( 1 - \frac{\delta_{\lambda_4} \lambda_4}{\log\log{\mathcal{H}}}\right) \frac{R_1 (\log{T})^{\lambda_3 - \lambda_1} (\log\log{T})^{\lambda_4 - \lambda_2}}{W \pi (\lambda_3 + 1)}
    + \frac{4 e^{c+1} \log{T}}{\log{x}} .
\end{split}
\end{equation}
\end{proposition}

\begin{proof}   
It follows from \eqref{eqn:approximation} and Lemma \ref{lem:important_integral_Chalker} that if $x > 1$, $c = 1 + 1/\log{x}$, and $T \geq \mathcal{H}$, then
\begin{equation}\label{eqn:not_asymptotic}
    | M(x) | 
    < \mathcal{E}_{R_1,R_2,\sigma_0}(x,T) + \left(2 e (\log{x} + 1) + 4 e^{c+1} \log{T} \right) \frac{x}{T} + 4e^c .
\end{equation} 
For each configuration of $\lambda_1$, $\lambda_2$, $\lambda_3$, and $\lambda_4$, we have $x^{\sigma_W(T)} = x/T$ and that $(\log{t})^{\lambda_1} (\log\log{t})^{\lambda_2}$, $(\log{t})^{\lambda_3} (\log\log{t})^{\lambda_3}$ are both increasing functions on $t \geq 14$. By definition, we also note that
\begin{align*}
    \frac{\mathcal{E}_{R_1,R_2,\sigma_0}(x,T)}{x/T}
    &= \frac{R_1}{\pi \log{x}} \left( \frac{e W (\log{T})^{\lambda_3 + 1} (\log\log{T})^{\lambda_4}}{W \log{T}} + \frac{x^{\sigma_W(\mathcal{H}) - 1} T (\log{\mathcal{H}})^{\lambda_3} (\log\log{\mathcal{H}})^{\lambda_4}}{\mathcal{H}} \right) \\
    &\qquad + \frac{R_1 x^{\sigma_W(T) - 1} T}{\pi} \Bigg( \left( 1 - \frac{\delta_{\lambda_4} \lambda_4}{\log\log{T}}\right) \frac{(\log{T})^{\lambda_3 + 1} (\log\log{T})^{\lambda_4}}{\lambda_3 + 1} \\ 
    &\hspace{8cm}- \frac{(\log{\mathcal{H}})^{\lambda_3 + 1} (\log\log{\mathcal{H}})^{\lambda_4}}{\lambda_3 + 1} \Bigg)  \\
    &\qquad + \frac{x^{\sigma_0 - 1} T}{\pi} \left( \frac{R_2 (\log{\mathcal{H}})^2}{2} + \frac{56}{\sigma_0(\sigma_0 - \tfrac{1}{2})} \right) .
\end{align*}
For our choice $T = T_x$ satisfying \eqref{eqn:T_choice}, we note $x^{\sigma_W(\mathcal{H}) - 1} T \leq x^{\sigma_W(T) - 1} T = 1$ and
\begin{align*}
    \frac{e W (\log{T})^{\lambda_3 + 1} (\log\log{T})^{\lambda_4}}{W \log{T} \log{x}} 
    = \frac{e (\log{T})^{\lambda_3 -\lambda_1} (\log\log{T})^{\lambda_4 - \lambda_2}}{W \log{T}} ,
\end{align*}
hence
\begin{align*}
    \frac{\mathcal{E}_{R_1,R_2,\sigma_0}(x,T)}{x/T}
    &\leq \frac{R_1}{\pi} \Bigg( \frac{e (\log{T})^{\lambda_3 -\lambda_1 - 1} (\log\log{T})^{\lambda_4 - \lambda_2}}{W} \\
    &\hspace{2cm} + \left( 1 - \frac{\delta_{\lambda_4} \lambda_4}{\log\log{\mathcal{H}}}\right) \frac{(\log{T})^{\lambda_3 -\lambda_1} (\log\log{T})^{\lambda_4 - \lambda_2} \log{x}}{W(\lambda_3 + 1)} \\ 
    &\hspace{2cm} + \frac{(\log{\mathcal{H}})^{\lambda_3} (\log\log{\mathcal{H}})^{\lambda_4}}{\mathcal{H} \log{x}} - \frac{(\log{\mathcal{H}})^{\lambda_3 + 1} (\log\log{\mathcal{H}})^{\lambda_4}}{\lambda_3 + 1} \Bigg) \\ 
    &\quad + \frac{\exp\left( \log{T} + (\sigma_0 - 1) \log{x}\right)}{\pi} \left( \frac{R_2 (\log{\mathcal{H}})^2}{2} + \frac{56}{\sigma_0(\sigma_0 - \tfrac{1}{2})} \right) .
\end{align*}
Note that $(\log{T})^{\lambda_3 - \lambda_1 - 1} (\log\log{T})^{\lambda_4 - \lambda_2}$ decreases on $T \geq \mathcal{H}$ in every case. 
Assert this upper bound in \eqref{eqn:not_asymptotic} and combine with the condition \eqref{eqn:suffices_to_show}, which holds for all $x \geq x_0$ by assumption, to recover the stated result.
\end{proof}

If $\lambda_1 = 1$, $\lambda_3 = 11/12$, and $\lambda_2 = \lambda_4 = 0$ in Proposition \ref{prop:MF_framework}, then the definition of $T$ in \eqref{eqn:T_choice} is explicit and can be directly applied. However, in all our other cases, the definition of $T$ is implicit, so we need bounds for $T$ depending on $x$ to apply Proposition \ref{prop:MF_framework} in a meaningful way. In fact, in Lemma \ref{lem:logT_bounds} below, we demonstrate the solution to the implicit definition of $T$ satisfies 
\begin{equation}\label{F def}
    \log{T_x} \asymp F_{W,\lambda_1,\lambda_2}(x) := \left(\frac{\log{x}}{W}\right)^{\frac{1}{\lambda_1 + 1}} \left(\frac{\lambda_1 + 1}{\log\log{x} - \log{W}}\right)^{\frac{\lambda_2}{\lambda_1 + 1}} .
\end{equation}

\begin{lemma}\label{lem:logT_bounds}
Assume the notation in Proposition \ref{prop:MF_framework} and fix admissible choices for $W$, $\lambda_1$, $\lambda_2$. We have
\begin{equation*}
    1 
    \leq \frac{\log{T}}{F_{W,\lambda_1,\lambda_2}(x)}
    \leq \left(1 + \frac{\lambda_2 \log\log{F_{W,\lambda_1,\lambda_2}(x)}}{(\lambda_1 + 1) \log{F_{W,\lambda_1,\lambda_2}(x)}}\right)^{\frac{\lambda_2}{\lambda_1 + 1}} .
\end{equation*}
\end{lemma}
\begin{proof}
We have asserted $T = T_{x}$ determined by \eqref{eqn:T_choice}, so together with \eqref{F def}, we have
\begin{align}\label{F/T}
    \frac{F_{W,\lambda_1,\lambda_2}(x)}{\log{T}} 
    = \left(1 + \frac{\lambda_2 \log\log\log{T}}{(\lambda_1 + 1) \log\log{T}}\right)^{-\frac{\lambda_2}{\lambda_1 + 1}} 
\end{align}
with $T \geq \mathcal{H}$. In particular, with $\lambda_2 < 0$ we observe
\begin{align*}
    \frac{F_{W,\lambda_1,\lambda_2}(x)}{\log{T}} 
    = \left(1 - o(1)\right)^{-\frac{\lambda_2}{\lambda_1 + 1}} ,
\end{align*} 
and with $\lambda_2 \geq 0$, we can replace $1 - o(1)$ with $1 + o(1)$. 
Therefore, treating the cases $\lambda_2 < 0$ and $\lambda_2 \geq 0$ independently, it follows that
\begin{align*}
    \log{T} &\geq F_{W,\lambda_1,\lambda_2}(x) 
\end{align*}
for $T\ge \mathcal{H}$ and all our pairs of $(\lambda_1,\lambda_2)$.

Similarly, treating the cases with $\lambda_2 < 0$ and $\lambda_2 \geq 0$ independently, we use the fact that  $\log\log\log{T} / \log\log{T}$ decreases for $T\ge \mathcal{H}$, and apply the above inequality to \eqref{F/T} to reveal the result. 
\end{proof}

\subsection{Proof of (\ref{eqn:classical-ob2})}\label{ssec:classical}

Assume the same notation as in Proposition \ref{prop:MF_framework}, fixing $$W \in \{5.56,6,7,8,10\},$$ $\lambda_1 = 1$, $\lambda_3 = 11/12$, $\lambda_2 = \lambda_4 = 0$. In addition, we let $\sigma_0 = 0.75$, $R_2 = e^{146}$, and $R_1 = \mathcal{R}_1$ defined for respective choices of $W$ as in Table \ref{tab:objectives}. These choices for $R_1$ and $R_2$ are admissible in light of Theorem \ref{thm:objective1} and Lemma \ref{lem:central_inequality_lwr}. 
For any $x \geq \exp(W(\log{\mathcal{H}})^{2})$, we have 
\begin{equation*}
    \log{T} = \sqrt{\frac{\log{x}}{W}} \geq \log{\mathcal{H}}
\end{equation*}
and observe that \eqref{eqn:suffices_to_show} holds. Therefore, the result \eqref{eqn:classical-ob2} follows upon asserting the above choices in Proposition \ref{prop:MF_framework}. \qed

\subsection{Proof of (\ref{eqn:Littlewood-ob2})}\label{ssec:Littlewood}

Assume the same notation as in Proposition \ref{prop:MF_framework}, fixing $$W \in \{21.24,21.5,22,25,30\},$$ $\lambda_1 = 1$, $\lambda_2 = -1$, $\lambda_3 = 11/12$, and $\lambda_4 = -3/4$. Since $T \geq \mathcal{H}$ and $W \geq 21.24$, it follows from \eqref{eqn:T_choice} that 
\begin{equation*}
    \log x 
    = \frac{W\log T}{\log\log T}\log T
    \geq \frac{21.24\log{\mathcal{H}}}{\log\log{\mathcal{H}}}\log T \ge \log T , 
\end{equation*}
or equivalently, $T \leq x$. Therefore, it subsequently follows from \eqref{eqn:T_choice}, $\lambda_1 = 1$, and $\lambda_2 = -1$ that
\begin{equation*}
    \log{T} 
    = \sqrt{\frac{\log x \log\log T}{W}}
    \leq \sqrt{\frac{\log x \log\log x}{W}} .
\end{equation*}
Next, let $\sigma_0 = 0.75$, $R_2 = e^{146}$, $R_1 = \mathcal{R}_2$, and $x_0 = x_2(W)$, where $\mathcal{R}_2$, $x_2(W)$ are defined for respective choices of $W$ as in Table \ref{tab:objectives}. 
These choices for $R_1$ and $R_2$ are admissible in light of Theorem \ref{thm:objective1} and Lemma \ref{lem:central_inequality_lwr}. To compute $x_2(W)$ for each fixed choice of $W$, we systematically incremented choices (rounded up to two decimal places) for $\log{x_2(W)}$ until $T \geq \mathcal{H}$ and \eqref{eqn:suffices_to_show} hold for every $x \geq x_2(W)$. 
The result \eqref{eqn:Littlewood-ob2} follows by substituting the above choices into Proposition \ref{prop:MF_framework}, and then combining this with the inequalities
\begin{equation*}
    \log{\mathcal{H}} \leq \log{T} \leq \left( \frac{\log{x} \log\log{x}}{W} \right)^{\frac{1}{\lambda_1 + 1}} ,
\end{equation*}
which are used to verify the sufficient condition \eqref{eqn:suffices_to_show}, together with Lemma \ref{lem:logT_bounds}. \qed

\subsection{Proof of (\ref{eqn:KV-ob2})}\label{ssec:KV}

Assume the same notation as in Proposition \ref{prop:MF_framework}, fixing $$W \in \{53.99,54,55,60,70\},$$ $\lambda_1 = \lambda_3 = 2/3$, $\lambda_2 = 1/3$, and $\lambda_4 = 1/4$. Since $T \geq \mathcal{H}$ and $W \geq 53.99$, it follows from \eqref{eqn:T_choice} that 
\begin{equation*}
    \log{T} 
    = (\log\log{T})^{-\frac{1}{5}} \left(\frac{\log{x}}{W}\right)^{\frac{3}{5}}
    < \left(\frac{\log{x}}{W}\right)^{\frac{3}{5}} .
\end{equation*}
Next, let $\sigma_0 = 0.75$, $R_2 = e^{146}$, $R_1 = \mathcal{R}_3$, and $x_0 = x_3(W)$, where $\mathcal{R}_3$, $x_3(W)$ are defined for respective choices of $W$ as in Table \ref{tab:objectives}. 
These choices for $R_1$ and $R_2$ are admissible in light of Theorem \ref{thm:objective1} and Lemma \ref{lem:central_inequality_lwr}. 
To compute $x_3(W)$ for each fixed choice of $W$, we systematically incremented choices (rounded up to two decimal places) for $\log\log{x_3(W)}$ until $T \geq \mathcal{H}$ and \eqref{eqn:suffices_to_show} hold for every $x \geq x_3(W)$. 
The result \eqref{eqn:KV-ob2} follows by substituting the above choices into Proposition \ref{prop:MF_framework}, and then combining this with the inequalities
\begin{equation*}
    \log{\mathcal{H}} \leq \log{T} < \left( \frac{\log{x}}{W} \right)^{\frac{3}{5}} ,
\end{equation*}
which are used to verify the sufficient condition \eqref{eqn:suffices_to_show}, together with Lemma \ref{lem:logT_bounds}. \qed

\section*{Acknowledgements}

We are grateful to Tim Trudgian for bringing the problem to our attention. In addition, we thank Dave Platt, Harald Helfgott, and other colleagues for their valuable feedback and insightful discussions. 
ESL thanks the Heilbronn Institute for Mathematical Research for their support. 
NL gratefully acknowledges that this research was supported in part by the Pacific Institute for the Mathematical Sciences.

\bibliographystyle{amsplain} 
\bibliography{references}

\end{document}